\documentclass{article}
\usepackage{amsfonts}
\usepackage{amsmath} 
\usepackage{amssymb}
\usepackage{graphicx}
\usepackage{amsthm}
\usepackage{mathtools}
\usepackage{mathrsfs}
\usepackage{hyperref}
\usepackage{dirtytalk}
\usepackage{textgreek}
\usepackage{wrapfig} 
\usepackage{tikz-cd}
\usepackage{shuffle}
\usepackage{xcolor}
\usepackage{euscript}
\usepackage{bbm}
\usepackage{booktabs}
\usepackage{enumitem}
\usepackage[a4paper, margin=1.2in]{geometry}

\usepackage{todonotes}
\setuptodonotes{backgroundcolor=white, textcolor=black, size=\tiny}

\usepackage{amsthm}
\usepackage{aliascnt}

\newtheorem{theorem}{Theorem}[section]

\newaliascnt{lemma}{theorem}
\newtheorem{lemma}[lemma]{Lemma}
\aliascntresetthe{lemma}

\newaliascnt{corollary}{theorem}
\newtheorem{corollary}[corollary]{Corollary}
\aliascntresetthe{corollary}

\newaliascnt{proposition}{theorem}
\newtheorem{proposition}[proposition]{Proposition}
\aliascntresetthe{proposition}

\newaliascnt{conjecture}{theorem}

\aliascntresetthe{conjecture}

\theoremstyle{definition}
\newaliascnt{definition}{theorem}
\newtheorem{definition}[definition]{Definition}
\aliascntresetthe{definition}

\newaliascnt{example}{theorem}

\aliascntresetthe{example}

\newaliascnt{notation}{theorem}
\newtheorem{notation}[notation]{Notation}
\aliascntresetthe{notation}

\theoremstyle{remark}
\newaliascnt{remark}{theorem}
\newtheorem{remark}[remark]{Remark}
\aliascntresetthe{remark}

\def\dif{\mathrm{d}}

\def\e{\mathrm{e}}
\def\b{v}
\def\c{w}

\date{\today}
\author{Emilio Ferrucci\thanks{Mathematical Institute, University of Oxford. Andrew Wiles Building, Woodstock Rd, Oxford OX2 6GG UK.} \and Timothy Herschell\footnotemark[1] \and Christian Litterer\thanks{Department of Mathematics, University of York. James College, York YO10 5DD, UK.} \and Terry Lyons\footnotemark[1]}
\title{High-degree cubature on Wiener space \\ through unshuffle expansions}

\begin{document}

\maketitle
\begin{abstract}
Utilising classical results on the structure of Hopf algebras, we develop a novel approach for the construction of cubature formulae on Wiener space based on unshuffle expansions. We demonstrate the effectiveness of this approach by  constructing the first degree-7 cubature formula on $d$-dimensional Wiener space with drift in the sense of Lyons and Victoir \cite{LV04} which is explicit as a function of an underlying Gaussian cubature. The support of our degree-7 formula is significantly smaller than that of currently implemented or proposed constructions.
\end{abstract}

\section*{Introduction}
\setcounter{section}{-1}
\refstepcounter{section}\label{sec:intro}

Cubature, in combination with Taylor expansion for error estimation, is a
classical and efficient method for approximating the integrals of
sufficiently regular functions of several variables. The cubature paradigm replaces a target
measure with a discrete measure of small, finite support that exactly
matches the integrals of a finite-dimensional space of (polynomial) test
functions. Cubature on Wiener space (as developed by Lyons and Victoir \cite%
{LV04}) extends this concept to path-space: the traditional Taylor
expansion is replaced by the stochastic Taylor expansion, and
Stratonovich iterated integrals of Brownian motion take the place of
polynomials. In its iterated form, known as the Kusuoka-Lyons-Victoir (KLV)
method (see Lyons-Victoir \cite{LV04} and Kusuoka \cite{kusuoka}),
cubature on Wiener space provides high-order weak approximations of measures
evolving in infinite-dimensional spaces.\\

The KLV method, in its basic form, approximates a broad class of potentially
hypo-elliptic parabolic partial differential equations (PDEs) by computing
weak solutions of stochastic differential equations (SDEs). Unlike classical
cubature methods, it can be shown using Malliavin calculus that these
approximations achieve high-order accuracy even when the test function only
has Lipschitz regularity. This foundational work has opened up a range of
cubature applications, including McKean-Vlasov SDEs \cite{McKeanVlasov}, \cite{CMM}, backward stochastic
differential equations \cite{CM1}, \cite{CM2}, sensitivity analysis for financial derivatives \cite{T},
stochastic partial differential equations \cite{BT}, and the non-linear filtering
problem \cite{CG}, \cite{CO}. Numerical implementations of cubature on Wiener space combining the
KLV\ method with partial sampling techniques in Gyurko, Lyons \cite{gyurko}
and recombination in Ninomiya, Shinozaki \cite{ninomiya} demonstrate the
high-order convergence properties of cubature up to
degree nine. These implementations also show that high-degree cubature
approximations can be significantly more efficient than first- and
second-order methods, especially when high-accuracy solutions are required.
A numerical toy example implemented in \cite{LL12} used high-degree Gaussian quadrature to make the recombining KLV method adaptive to
the regularity of the boundary data. Although this approach produces
exceptionally accurate and efficient approximations of the PDE solution, the
high-degree cubature measures on path-space necessary for extending this
methodology to practical problems remain unavailable.
\\

The construction of efficient, high-degree cubature measures for more than
two- or three-dimensional noise remains a significant challenge for Wiener
space. The algebraic and computational complexity increases rapidly with the
degree of the approximation, see e.g.\ \cite{gyurko}. 
\\

Explicit constructions of cubature measures on Wiener space broadly involve
three steps.
\begin{enumerate}
\item Pick a convenient basis of the free Lie algebra and expand the expected signature of Brownian motion (augmented with drift) as a symmetric tensor product of these basis elements.
\item Define a suitable Lie polynomial with unknown, random coefficients. Setting its expected tensor-exponential equal to the expected signature of desired degree gives rise to a moment problem.
\item Finally, the most challenging step. Solve this moment problem by constructing the random coefficients as polynomials of (a small family of) distributions for which (efficient) finite-dimensional joint cubature measures of the necessary orders are available. Replacing these random variables with their cubature measures yields the cubature measure on Wiener space.
\end{enumerate}
Crucially, producing a cubature formula does not only entail solving a moment problem but also \emph{designing} one (with as few degrees of freedom as possible) so that it \emph{can} be solved. For this reason, steps 2 and 3 must often be iterated many times before a good solution is found.

In this paper we propose a novel approach based on unshuffle expansions that leverages the Hopf algebra structure of the tensor algebra. Instead of a basis of the free Lie algebra, we use a redundant spanning set given by the images of the linear
extension of the logarithm restricted to grouplike elements, in order to expand the expected
signature. This map is also known as the canonical projection onto the free
Lie algebra \cite{Reu03}, or Eulerian idempotent \cite{Lod94,Pat93}, and occurs
naturally in the non-commutative correction terms of the
Campbell-Baker-Hausdorff formula. The choice to abandon the use of a basis, while counterintuitive, has the benefit of leveraging the symmetries of the Wiener measure in a way that makes the construction more tractable.

Our approach has several advantages over existing methods in \cite{shinozaki,hayakawa-tanaka,litterer, gyurko}. The unshuffle expansion of the expected signature in this spanning set leads
to a simplified and sparser moment problem. This greatly mitigates the technical complexity that makes existing approaches increasingly
intractable for higher degrees. As a consequence, we can construct a degree-seven cubature measure on Wiener space that, in contrast to existing constructions, can be stated explicitly for arbitrary dimensions in terms of an underlying Gaussian cubature of inhomogeneous degree. This results in measures with vastly smaller and, hence, more efficient supports than existing constructions. 

Our paper is organised as follows. In \autoref{sec:background} we review background notions: namely, \autoref{subsec:prelims} is a gentle introduction to some well-known facts about the tensor Hopf algebra, and \autoref{subsec:cubpoly} is a reminder on cubature on Wiener space using the formulation in terms of Lie polynomials. Our original contribution is contained in \autoref{sec:generalconstruction}, which builds up to the main result \autoref{theorem:deg_7_cubature_formula}, a degree-7 formula on Wiener space in arbitrary dimension plus drift. This theorem is later refined in \autoref{cor:sym_rem}, and we end the section with a detailed comparison to existing constructions, \autoref{sec:gaussiancubatures}. The appendices \autoref{app:lemma_proofs}, \autoref{app:lyndon} contain the more technical aspects of the calculations, and \autoref{app:numerics} contains a numerical toy example showcasing the effectiveness of our formula. Before we begin, we provide some context for our results by reviewing the literature on existing methods for cubature on Wiener space.

\subsection*{Existing constructions}

Constructions for cubature measures on Wiener space of degree three and five
were first obtained in Lyons, Victoir \cite{LV04} by solving moment
problems for discrete random variables arising from the expansion of the
non-commutative exponential in a Poincare-Birkhoff-Witt basis of the free
tensor algebra. Since then, several attempts have been made to develop
constructions of cubature measures on Wiener space that extend to higher
degrees. For the one- and two-dimensional noise Litterer \cite{litterer}\
generalises the Lyons-Victoir construction to degree seven. Using a similar approach Gyurko and Lyons \cite{gyurko}
obtained measures of degree nine and eleven for one-dimensional noise.
{ In three dimensions, a degree-seven generalisation of the Lyons–Victoir cubature in the Hall basis has been obtained by Herschell \cite{timothy}, extending the original construction of \cite{LV04}. As observed by Gyurko–Lyons \cite{gyurko}, such a 3D degree-7 construction has support on 91 Lie-basis elements whose coefficients must satisfy more than 150 inhomogeneous polynomial constraints, underscoring the difficulty of the problem. The size of the support of the formula in \cite{timothy} serves as a key benchmark for our results; in particular, we recover the corresponding support size for three-dimensional noise from our construction in \autoref{sec:brokensym}. While
these constructions result in efficient, explicit cubature formulae on
Wiener space with the smallest support known to date they do not extend in
any obvious and tractable way to higher dimensions due to the lack of symmetry in both the basis and the associated moment constraints. }

The first construction of a general, if not explicit, degree-seven
cubature measure on Wiener space is due to Shinozaki \cite[Theorems 3.1 and 3.4]{shinozaki}, who leverages the algebraic relations of
products of iterated stochastic integrals to construct a moment similar
family on grouplike elements that matches the expectations of Stratonovich
iterated integrals of Brownian motion up to degree seven. The construction
is a remarkable technical achievement, but does not lead to a cubature
measure on Lie polynomials or paths that can be explicitly stated in
dimensions greater than two. \cite[Theorem 3.4]{shinozaki} is stated for the two-dimensional
case and the intrinsic technical complexity of the calculations means it
cannot be stated in higher dimensions. Despite this, higher-dimensional
examples can at least in principle be obtained from Theorem 3.1 using
computer algebra programs \cite[Remark 3.7]{shinozaki}. Shinozaki's construction
relies on high-dimensional Gaussian cubature measures. The $7$-moment similar
families for $d$-dimensional Brownian motion requires discrete random variables matching sufficient moments of a 
\begin{equation*}
2d+3 {d \choose 2} +3{ d \choose 3} =\frac{d^{3}+3d }{2}, \quad d\geq 3
\end{equation*}%
dimensional standard normal random variable. For the case of two-dimensional Brownian motion,
Shinozaki's construction has been implemented using a seven-dimensional
degree-seven Gaussian cubature (Ninomiya, Shinozaki \cite{ninomiya}).

An alternative, randomised construction of cubature measures based on
recombination was first proposed in Litterer, Lyons \cite[p.1307-1308]{LL12}. The algorithm is based on an observation by Wendel \cite{wendel} who observed that for a spherically symmetric measure, an IID sample 
of $2n$ points contains the origin inside its convex hull with probability $1/2$. In Litterer \cite{litterer} the signature was also used as a
test function for recombination on path-space. A sophisticated extension of
these ideas in Hayakawa, Tanaka \cite{hayakawa-tanaka} and Hayakawa, Oberhauser, Lyons \cite{hayakawa-lyons} provides far more general estimates
that also apply to the number of paths required for the randomised
construction of cubature measures on Wiener space. Since the randomised
construction makes no use of the symmetries inherent in the Wiener measure,
the support of the resulting formula grows with the dimension of the
truncated tensor algebra over $d$ variables. The computational
complexity in the examples considered grows even faster, making the proposed
construction for degree-seven and above numerically intractable even for
moderate-dimensional Brownian noise. For degree seven, only a two-dimensional example is constructed in 
\cite{hayakawa-tanaka}.

The approach proposed in this paper resolves some significant limitations of
existing constructions. It allows us to construct explicit cubature measures
for $d$-dimensional Brownian motion based on a $d$-dimensional degree-$7$
Gaussian cubature measures and some auxiliary variables that only have to
match Gaussian moments up to degree three. This results in cubature measures
on Lie polynomials that are explicit for arbitrary dimension and have in many important cases vastly
smaller support compared to any existing
construction (cf.\ \autoref{sec:effic} for a detailed discussion). 

The support of any cubature measure
on Wiener space is trivially bounded below by the size of the support of its
underlying $d$-dimensional Gaussian cubature measure (cf.\ Lyons, Victoir \cite%
{LV04}). We can show that the size of the support of our cubature measure for $d$-dimensional noise is bounded above by $4d^2$ times the size of the support of any $d$-dimensional degree-seven Gaussian cubature measure (cf.\ \autoref{sec:gaussiancubatures} for a discussion of such measures).

\section{Background on algebra and cubature}\label{sec:background}

\subsection{The tensor Hopf algebra and its Eulerian idempotent}\label{subsec:prelims}

In this section we provide the algebra background which will be used for the construction of the cubature formula; for details we refer to \cite{Reu03}.

Let $V$ be a finite-dimensional vector space. We denote the tensor algebra
\[
T(V) \coloneqq \bigoplus_{n = 0}^\infty V^{\otimes n},
\]
endowed with the tensor product $\otimes$, making $T(V)$ the free (non-commutative) algebra generated by $V$. As a vector space, it is spanned by elementary tensors $v_1 \otimes \cdots \otimes v_n$ with $v_k \in V$, which we identify with \emph{words} $v_1\ldots v_n$, omitting the tensor product symbol; we use $1$ to denote the empty word, the generator of $\mathbb R = V^{\otimes 0}$, and will sometimes call elements of $V$ \emph{letters}. The \emph{unshuffle coproduct} is defined by $\Delta_\shuffle 1 \coloneqq 1 \otimes 1$ and
\begin{equation}\label{eq:coproduct}
\begin{split}
\Delta_\shuffle \colon T(V) &\to T(V) \otimes T(V) \\ 
v_1\ldots v_n &\mapsto \sum_{I \sqcup J = [n]} v_I \otimes v_J, \qquad v_1,\ldots,v_n \in V,
\end{split}
\end{equation}
where we are summing over partitions of the set with $n$ elements into two sets $I$ and $J$, and $v_K \coloneqq v_{k_1}\ldots v_{k_p}$ for $K = \{k_1,\ldots,k_p\}$ with $k_1 < \ldots < k_p$. In other words, $\Delta_\shuffle$ separates a word $v_1\ldots v_n$ into two subwords without altering the order they inherit. For example,
\begin{equation}\label{ex:abc}
\Delta_\shuffle (u\b\c) = 1 \otimes u \b \c  + u \otimes \b \c  + \b  \otimes u \c  + \c  \otimes u \b  + \b \c  \otimes u + u \c  \otimes \b  + u \b  \otimes \c  + u \b \c  \otimes 1 .
\end{equation}
The number of terms in the expression for the coproduct of a word grows very rapidly; for example $\Delta_\shuffle (u\b\c z)$ is a sum of 16 terms which include ones such as $\b \otimes u\c z$ and $uz \otimes \b\c$. The coproduct is coassociative, i.e.\ its iterates are well-defined independently of the order: define
\[
\Delta_\shuffle^1 \coloneqq \mathbbm 1_{T(V)}, \quad \Delta^n_\shuffle \coloneqq (\mathbbm 1 \otimes \Delta^{n-1}_\shuffle) \otimes \Delta_\shuffle = ( \Delta^{n-1}_\shuffle \otimes \mathbbm 1) \otimes \Delta_\shuffle \colon T(V) \to T(V)^{\otimes n},
\]
for $n \geq 2$ (so that in particular $\Delta_\shuffle = \Delta_\shuffle^2$). It can be shown that $\Delta_\shuffle$ is also an algebra morphism, and that $\otimes$ is a coalgebra morphism, making $(T(V), \otimes, \Delta_\shuffle)$ a bialgebra {\cite[Proposition 1.9]{Reu03}}; since it is graded and connected, it can be endowed with a (unique) antipode, making it a Hopf algebra, the \emph{tensor Hopf algebra}. It is the dual Hopf algebra to the (perhaps better known) \emph{shuffle Hopf algebra} $(T(V), \shuffle, \Delta_\otimes)$ in which the coproduct is defined by deconcatenation, i.e.\ splitting a word in all possible ways. One should think of the shuffle Hopf algebra as indexing iterated integrals (with $\shuffle$ encoding the operation of product of two iterated integrals), while the tensor Hopf algebra as indexing vector fields (with the tensor product encoding composition of vector fields). We will always be working in the tensor Hopf algebra.

\begin{remark}\label{rem:cocommutative}
It is useful to remark that $\Delta_\shuffle$ is cocommutative, i.e.\ $\tau \circ \Delta_\shuffle = \Delta_\shuffle$ where $\tau$ swaps the two copies of $T(V)$. Therefore, it is possible to express $\Delta_\shuffle^n$ by summing over ordered submultisets and using symmetric products, see \eqref{eq:symmprod} below, instead of subwords (for example, in \eqref{ex:abc}, $v \otimes uw$ and $uw \otimes v$ can be grouped together). This is advantageous in implementations, as it reduces the number of terms factorially.
\end{remark}

{Notice that $\Delta_\shuffle x$ is a sum which always contains the two summands $1 \otimes x$ and $x \otimes 1$; a similar comment also applies to higher-order coproducts, in which some summands have the degree-0 element $1 \in \mathbb R = V^{\otimes 0}$ in one of the \say{slots}. Sometimes it is helpful to remove these trivial summands, and for this purpose one therefore defines} the (iterated) \emph{reduced} coproduct by
\[
\widetilde \Delta_\shuffle^n \coloneqq \pi_{\geq 1}^{\otimes n} \circ \Delta_\shuffle^n, \qquad \text{where } \pi_{\geq 1} \colon T(V) \twoheadrightarrow \bigoplus_{n = 1}^\infty V^{\otimes n}
\]
is the projection onto tensors of positive degree (for example $\widetilde \Delta_\shuffle x = \Delta_\shuffle x - x \otimes 1 - 1 \otimes x$). We introduce similar notation for the projections of the tensor algebra onto its graded components
\begin{equation*}
\pi_n \colon T(V) \twoheadrightarrow \bigoplus_{k = 0}^n V^{\otimes k}.
\end{equation*}
It will be helpful to use sum-free Sweedler notation
\[
\Delta_\shuffle^n { x} \eqqcolon  x_{(1)} \ldots x_{(n)}, \quad \widetilde\Delta_\shuffle^n x \eqqcolon  x^{(1)} \ldots x^{(n)} .
\]
To be more precise, anytime $x_{(1)}, \ldots, x_{(n)}$ appear in some expression, this means we are taking the $m$-fold coproduct of $x$ and performing some operation on the individual factors (and similarly in the reduced case, with superscripts instead of subscripts); we will provide an example of this notation shortly \eqref{eq:star}. For $x_1,\ldots,x_n \in T(V)$ define the symmetric (tensor) product by 
\begin{equation}\label{eq:symmprod}
(x_1 , \ldots ,  x_n ) \coloneqq \frac{1}{n!} \sum_{\sigma \in \mathfrak S_n} x_{\sigma(1)} \otimes \cdots \otimes x_{\sigma(n)} \ \in T(V).
\end{equation}
Note that the symmetric product is not associative, e.g.\ the above product cannot be computed recursively by $((x_1 , \ldots , x_{n-1}) , x_n)$.

The \emph{convolution product} on linear endomorphisms of the vector space $T(V)$ is defined by
\[
\star \colon \mathrm{End}(T(V))^{\otimes 2} \to \mathrm{End}(T(V)), \qquad f \star g \coloneqq \otimes \circ (f \otimes g) \circ \Delta_\shuffle .
\]
Here \say{$\otimes \circ$} denotes the product mapping $T(V)^{\otimes 2} \to T(V)$, and note that by cocommutativity this can be replaced by the symmetric tensor product. $\mathrm{End}(T(V))$ forms a group under $\star$, with neutral element the projection $\pi_0 \colon T(V) \twoheadrightarrow V^{\otimes 0}$. Thanks to the bialgebra properties, multiple convolution products can be written using the iterated coproduct:
\begin{equation}\label{eq:star}
\begin{split}
f_1 \star \cdots \star f_n &= \otimes^n \circ (f_1 \otimes \cdots \otimes f_n) \circ \Delta_\shuffle^n, \\
\text{in Sweedler notation:}\quad (f_1 \star \cdots \star f_n)(x) &= f_1(x_{(1)}) \otimes \cdots \otimes f_n(x_{(n)}) \\
&= (f_1(x_{(1)}), \cdots , f_n(x_{(n)})).
\end{split}
\end{equation}
For $f \in \mathrm{End}(V)$, consider its exponential and logarithmic series under the convolution product:
\begin{align*}
\exp_\star(f) \coloneqq \sum_{n \geq 0} \frac{f^{\star n}}{n!}, \qquad \log_\star (f) \coloneqq \sum_{n \geq 1} \frac{(-1)^{n-1}}{n}(f - \pi_0)^{\star n} .
\end{align*}
When applied to any given element, these series always reduce to finite sums by conilpotency (the property that for any $x \in T(V)$ there exists $n$ such that $\Delta^n_\shuffle x = 0$).
\begin{definition}[Eulerian idempotent]\label{def:eul}
The linear map
\[
\e \coloneqq \log_\star(\mathbbm 1_{T(V)}) \in \mathrm{End}(V)
\]
is called the \emph{Eulerian idempotent} of the Hopf algebra $(T(V), \otimes, \Delta_\shuffle)$.
\end{definition}
The Eulerian idempotent is denoted $\pi_1$ in \cite[\S 3.2]{Reu03} (not to be confused with the map $\pi_1$ as it is denoted here, the projection onto single letters); we also refer to \cite{Pat93, Lod94} for further details.

Define a Lie bracket on $T(V)$ by $[x,y] \coloneqq x \otimes y - y \otimes x$. The free Lie algebra $\mathcal L(V)$ over $V$ is the smallest Lie subalgebra of $T(V)$ which contains $V$; it can be explicitly described as the direct sum
\[
\mathcal L(V) = V \oplus [V,V] \oplus [[V,V],V] \oplus \ldots,
\]
where Lie brackets between vector subspaces of $T(V)$ denote the space spanned by all Lie brackets of the two vector spaces; it can also be described more abstractly by a universal property, without reference to $T(V)$.
It is a non-trivial result that $\mathcal L(V)$ coincides with the space of $\Delta_\shuffle$-primitive elements {\cite[Theorem 3.1]{Reu03}}: 
\[
\mathcal L(V) = \{x \in T(V) \mid \widetilde \Delta_\shuffle x = 0\} .
\]
The expression for the Eulerian idempotent can be made more explicit by
\begin{align*}
\e(x) &= \sum_{n \geq 1} \frac{(-1)^{n-1}}{n} (\mathbbm 1 - \pi_0)^{\star n}(x) \\
&= \sum_{n \geq 1} \frac{(-1)^{n-1}}{n} \otimes^n \circ (\mathbbm 1 - \pi_0)^{\otimes n} \circ \Delta_\shuffle^n(x) \\
&= \sum_{n \geq 1} \frac{(-1)^{n-1}}{n} ( x^{(1)} , \ldots , x^{(n)}) .
\end{align*}
In the first identity we have used the bialgebra property to express convolution powers, and in the second we have used cocommutativity to symmetrise the tensor product (and note that the Sweedler notation is reduced).

We mention for the sake of completeness, though it will not be used in the following pages, that $x \in T(\!(V)\!) \coloneqq \prod_{n = 0}^\infty V^{\otimes n}$ is \emph{grouplike} if $\Delta_\shuffle x = x \otimes x$, and denote the group of these $\mathcal G(\!(V)\!)$. Calling $\mathcal G^n(V) \coloneqq \pi_n(\mathcal G(\!(V)\!))$, it still holds that for $x \in \mathcal G^n(V)$, $\Delta_\shuffle x = x \otimes x$ on $T^n(V)$. Therefore, calling $\mathcal G(V) \coloneqq \bigcup_{n \geq 0} \mathcal G^n(V) \subset T(V)$, we have that for $x \in \mathcal G(V)$
\begin{align*}
\e(x) = \sum_{n \geq 1} \frac{(-1)^{n-1}}{n} \otimes^n \circ (\mathbbm 1 - \pi_0)^{\otimes n} \circ \Delta_\shuffle^n(x) = \sum_{n \geq 1} \frac{(-1)^{n-1}}{n} x^{\otimes n} = \log (x)\,,
\end{align*}
i.e.\ $\e|_{\mathcal G(V)} = \log|_{\mathcal G(V)}$, where the latter denotes the logarithmic series taken w.r.t.\ the tensor product. In fact, $\e$ is the unique linear map $T(V) \to T(V)$ with this property, where uniqueness follows from the fact that $\mathcal G(V)$ spans $T(V)$. Recall that the inverse of $\log$ is the tensor exponential $\exp$, whose expression on a finite sum can be written in terms of the symmetric product as follows, cf.\ \cite[Proposition 4.3]{LV04}. For $x_1,\ldots,x_n \in T(V)$:
\begin{equation}\label{eq:exp}
\begin{split}
\exp \Big( \sum_{k = 1}^n x_k \Big) &= \sum_{m \geq 0} \sum_{k_1,\ldots, k_m = 1}^n \frac{1}{m!} (x_{k_1}, \ldots, x_{k_m}) \\ &= \sum_{\substack{m \geq 0 \\ m_1 + \ldots + m_n = m}} \frac{(\overbrace{x_1,\ldots, x_1}^{m_1}, \ldots, \overbrace{x_n,\ldots,x_n}^{m_n})}{m_1! \cdots m_n!}.
\end{split}
\end{equation}

We summarise the main properties of $\e$ of relevance to us in the following proposition{, which is classical (cf.\ \cite[Proof of Theorem 3.7]{Reu03}), and closely related to the celebrated Milnor-Moore theorem \cite{MM65}.}

\begin{proposition}\label{proposition:eul_id_projection}
The Eulerian idempotent is a projection onto $\mathcal L(V)$, and $x \in T(V)$ can be expressed as a symmetric product of Lie elements by
\begin{equation}\label{eq:expansion}
x = \sum_{n \geq 0} \frac{1}{n!} {(\e(x^{(1)}),\dots ,\e(x^{(n)}))} .
\end{equation}
In other words, this provides an isomorphism
\begin{equation}\label{eq:iso}
\bigoplus_{n = 0}^\infty \frac{\e^{\star n}}{n!} \colon T(V) \xrightarrow{\cong} \bigoplus_{n = 0}^\infty \mathcal L(V)^{\odot n}.
\end{equation}
\begin{proof}
Writing $\mathbbm 1 = \exp_\star (\log_\star(\mathbbm 1))$ yields the identity
\[
x = \exp_\star (\e)(x) = \sum_{n \geq 0} \frac{\e^{\star n}(x)}{n!} ,
\]
which coincides with the expression in the first statement, thanks to \eqref{eq:star}.
\end{proof}
\end{proposition}
We remark that the higher convolution powers of $\e$ that appear in \eqref{eq:expansion} have a simplified expression in terms of Stirling numbers of the first kind \cite{zbMATH00059889} (see \cite[Theorem 4.1.1]{CP21} for a recent presentation). We will also use the following two symmetries of the Eulerian idempotent. The first is a direct consequence of the isomorphism \eqref{eq:iso}.

\begin{proposition}[The symmetric property]\label{prop:symmetric}Let $x \in \bigoplus_{n = 2}^\infty V^{\odot n}$. Then $\e(x) = 0$. 
\end{proposition}
\noindent For the second symmetry, we introduce the \emph{reversal} operator, defined on elementary tensors as
\begin{equation}\label{eq:reversal}
* \colon T(V) \to T(V), \qquad (v_1 \ldots v_n)^* \coloneqq v_n \ldots v_1 \quad \text{for } v_1,\ldots, v_n \in V ,
\end{equation}
and extended linearly; $1^* = 1$.

\begin{proposition}[The reversal property, {\cite[Proposition 20]{burgunder}}]\label{prop:reversal} For $n \geq 1$ and $x \in V^{\otimes n}$
\[
\e(x) = (-1)^{n-1} \e(x^*) .
\]
\end{proposition}

Recall that a \emph{Hall basis} is a particular type of basis of the $\mathcal L(V)$. The following is an immediate consequence.
\begin{corollary}
Given a Hall basis $H$, the set
\begin{equation}\label{eq:PBW}
\{(h_1,\ldots,h_n) \mid n \in \mathbb N, \ h_1,\ldots,h_n \in H\}
\end{equation}
is a basis of $T(V)$, called the symmetrised Poincaré-Birkhoff-Witt (PBW) basis.
\end{corollary}

\begin{remark}
We take a moment to reflect on the difference between \eqref{eq:expansion} and the expression of $x$ in a basis \eqref{eq:PBW}. The latter is the choice made in \cite{LV04} and has the benefit that the expression of $x \in T(V)$ as a commutative polynomial in the Hall elements is unique. On the other hand, there is no canonical choice of a Hall basis, and such bases are inherently designed to break the symmetry; for example, the Lyndon Hall basis depends on the choice of a frame on $V$. As a consequence, expansions in Hall bases generally involve little structure and are therefore not easy to present concisely. On the other hand, \eqref{eq:expansion}, while not an expansion in a particular basis of $\mathcal L(V)$---the set $\{\e(w) \mid w \text{ word}\}$ is necessarily a dependent spanning set of $\mathcal L(V)$, even if the letters in $w$ are only taken to belong to a basis of $V$---has the property of only using the intrinsic structure of the Hopf algebra $(T(V), \otimes, \Delta_\shuffle)$. One significant benefit is that high-degree identities become tractable and presentable in a fraction of the space, thanks to the use of compact algebraic notation. Moreover, any formula that is mathematically derived using \eqref{eq:expansion} can always be implemented in a PBW basis. To illustrate this point, in \autoref{app:lyndon} we compute the images of the Eulerian idempotents of words of degree up to 5 (those needed in \autoref{theorem:deg_7_cubature_formula}) in the Lyndon Basis.
\end{remark}

Throughout this paper, $V$ will be $\mathbb R^{1+d}$, and we denote $\epsilon_0,\epsilon_1,\ldots,\epsilon_d$ its canonical frame. The zero-th coordinate will be reserved for the drift, which has twice the regularity of Brownian motion, and will therefore be given double the weight as the other coordinates, thus resulting in an inhomogeneously graded tensor algebra. It will be convenient to define
\begin{equation}\label{eq:pin-tilda}
\widetilde{\pi}_m \colon T(\mathbb R^{1+d}) \twoheadrightarrow \bigoplus_{2i + j \leq m} (\mathbb R^{1+d})^{\otimes (i,j)} ,
\end{equation}
where $(\mathbb R^{1+d})^{\otimes (i,j)}$ denotes the space spanned by all words of length $i + j$ of which exactly $i$ letters are the letter $0$. The inhomogeneous grading of the tensor algebra reflects the fact that the drift scales differently from Brownian motion. Because of the central role that the coordinates $\{\e(w) \mid w \text{ word}\}$ will play in this paper, we will use the following shorthand.

\begin{notation}\label{not:xi}
For $i_1,\ldots,i_n \in \{1,\ldots,d\}$ we denote
\begin{equation}
\xi_{i_1\ldots i_n} \coloneqq \e(\epsilon_{i_1} \ldots \epsilon_{i_n}).
\end{equation}
\end{notation}

\subsection{Cubature measures for Wiener space supported on Lie polynomials}\label{subsec:cubpoly}

Throughout this paper, $B$ will denote a $d$-dimensional Wiener process and $\widehat B$ is $B$ augmented with time in its zero-th coordinate, $\widehat B_t = (t, B_t)$. We write $S(\circ B)$ and $S(\circ \widehat B)$ for their respective Stratonovich signatures. A cubature measure (or formula) of degree $m$ for a probability measure $\rho$ on $\mathbb{R}^d$ is a finitely supported positive measure whose moments up to order $m$ agree with those of $\rho$. We recall the following definition of Lyons and Victoir for cubature measures on path space.

\begin{definition}[\cite{LV04}{Definition 2.2}]
\label{cubsigdef} We say a discrete probability measure $Q=\sum_{j=1}^{n}%
\lambda _{j}\delta _{\omega _{j}}$ supported on $n$
paths $\omega_j \in C^{1\text{-var}}([0,1], \mathbb R^{1+d})$ is a \emph{degree-$m$ cubature measure on $d$-dimensional Wiener space with drift}, if 
\begin{equation}
\mathbb{E}[\widetilde\pi_{m}(S_{0,1}(\circ \widehat B)) ]=\sum_{j=1}^{n}\lambda _{j}\widetilde\pi
_{m}(S_{0,1}(\omega _{j})),  \label{cubature}
\end{equation}
where $\widetilde\pi_m$ is defined in \eqref{eq:pin-tilda}.
\end{definition}
By scaling and stationarity of increments, a cubature measure on the interval $[s,t]$ can be obtained from $Q$ above by letting $\omega_{s,t;i}^{j}(u)=\sqrt{t-s}\cdot \omega
_{i}^{j}(u/(t-s))$, $j=1,\ldots ,d$ and keeping the weights of $Q$. The drift component of the cubature paths is given by $\omega^0(t) = t$. 

In the KLV method, cubature measures are used to weakly approximate solutions
to $\mathbb R^e$-valued Stratonovich stochastic differential equation (SDE) of the form 
\begin{equation}\label{eq:SDE}
\dif X_{t}=F(X_t) \circ \dif \widehat B_{t},
\end{equation}
with sufficiently regular vector fields $F$ which contains drift $F_0$ and the diffusion coefficients $F_k$, $k = 1,\ldots, d$. The KLV approximation is
computed by solving differential equations controlled by the cubature paths $%
\omega _{j}$ (rescaled according to the interval on which the equation is defined)
\begin{equation}\label{eq:CDE}
\dif Y_{j;t}=F(Y_{j;t}) \dif \omega_j(t) .
\end{equation}
For suitable test functions $f$, it is then shown \cite[Proposition 3.2]{LV04} that on the interval $[s,t]$
\[
\Big| \mathbb E[f(X_t) | X_s = x] - \sum_{j = 1}^n \lambda_j f(Y_{j;t}) \Big| \lesssim (t-s)^{(m+1)/2}, \qquad \text{with } Y_{j;s} = x, 
\]
with the constant of proportionality independent of $t-s$. Thanks to the independence of Brownian increments, this one-step estimate can be iterated along a partition to obtain a numerical scheme converging at rate $k^{1- (m+1)/2}$, where $k$ is the number of intervals in the partition \cite[Theorem 3.3]{LV04}.

It is well-known (e.g.\ \cite{kunita, benArous, castell_gaines, kusuoka}) that the one-step order-$m$ Taylor expansion of \eqref{eq:CDE}
on $[s,t]$ is equal, up to an error of order $(t-s)^{(m+1)/2}$, to the ODE
\begin{equation}\label{eq:log-ODE}
\dot Z_t =  F^\circ(Z) \widetilde \pi_m(\log S_{s,t}(\omega_j)), \quad Z_s = Y_{j;s} ,
\end{equation}
where the map $F^\circ:\mathcal L(\mathbb R^{1+d}) \to C^\infty(\mathbb R^e, \mathbb R^e)$ is the {restriction to $\mathcal L(\mathbb R^{1+d}) \subset T(\mathbb R^{1+d})$ of} the universal extension of the linear map that takes basis elements $\epsilon_i \in \mathbb R^{1+d}$ to the vector fields $F_i$. In the expression $F^\circ(z) \ell$, $z$ is the argument in $\mathbb R^e$ and $\ell$ is the argument in $\mathcal L(\mathbb R^{1+d})$ on which $F^\circ$ acts linearly, so that $F^\circ(z)\ell \in \mathbb R^e$; the inverse order of the arguments is motivated by the fact that we think of first evaluating $F^\circ$ at a state $z$ and then contracting it with a Lie element. When defined on a word, $F^\circ$ equals the differential operator
\[
F^\circ(z)\epsilon_{i_1} \ldots \epsilon_{i_n} \coloneqq (\varphi \mapsto F_{i_1} \circ \cdots \circ F_{i_n}\varphi(z)) ,
\]
where $\circ$ denotes composition of vector fields, i.e.\ $F_{i_1} \circ \cdots \circ F_{i_n}\varphi(z)$ is recursively equal to $F_{i_1}$ acting on the function $y \mapsto F_{i_2} \circ \cdots \circ F_{i_n}\varphi(y)$ evaluated at the point $z$. When the algebra homomorphism $F^\circ$ is restricted to $\mathcal L(\mathbb R^{1+d})$ it takes values in the Lie algebra of vector fields on $\mathbb R^e$,  $C^\infty(\mathbb R^e, \mathbb R^e)$. The error incurred by replacing ODEs controlled by paths $\omega_j$ with autonomous ODEs defined by $F^\circ (\widetilde \pi_m(\log S(\omega_j)))$ matches the order of error already present in the cubature approximation. We naturally arrive at the following alternative definition of a cubature measure.

\begin{definition}[\cite{LV04}{Definition 4.9}]
\label{cubfreeliedef} Letting $\ell_j \coloneqq \widetilde \pi_m(\log S(\omega_j))$, we say that the discrete probability $Q=\sum_{j=1}^{n}\lambda
_{j}\delta_{\ell_{j}}$ measure on $\mathcal L(\mathbb R^{1+d})$ is a \emph{degree-$m$ cubature measure on $d$-dimensional Wiener space with drift} if
\begin{equation*}
\mathbb{E}[\widetilde \pi_{m}( S_{0,1}(\circ \widehat B))] = \mathbb E_{Q}[\widetilde\pi_{m}\exp(\ell)].
\end{equation*}
\end{definition}

Both definitions of cubature on Wiener space, when applied iteratively in the KLV method, have the same order of convergence, see e.g.\ Cass-Litterer \cite{litterer-cass}. The KLV method based on measures on Lie polynomials corresponds to a version of Kusuoka's algorithm, Kusuoka \cite{kusuoka}. The logarithmic signature maps cubature paths to Lie polynomials. Conversely, Chow's theorem guarantees the existence of continuous bounded variation paths with logarithmic signature matching any Lie polynomial (cf.\ Lyons, Victoir \cite{LV04}, where measures on Lie polynomials serve as a crucial intermediate step for the construction of cubature measures on paths). In the following we will adopt \autoref{cubfreeliedef}, which has been preferred in implementations of high order cubature measures \cite{gyurko,ninomiya}. Also, it has the advantage of being algebraic, and of already providing the building blocks for most numerical ODE solvers, the simplest of which is a Taylor scheme.

\section{Explicit, general-\texorpdfstring{$d$}{d}, degree-\texorpdfstring{$7$}{7} cubature measures through unshuffle expansions}\label{sec:generalconstruction}
We propose the following variation and refinement of the three-point plan explained in \autoref{sec:intro} for constructing a degree-7 cubature measure on dimension-\(d\) Wiener space with drift in the sense of \autoref{cubfreeliedef}.
\begin{enumerate}
\item Expand $\mathbb{E}[\widetilde \pi_{m}( S_{0,1}(\circ \widehat B))]$ using the Eulerian idempotent expansion \eqref{eq:expansion}.
\item Write an unknown $\mathcal L(\mathbb R^{1+d})$-valued random variable $\mathcal L$ as a linear combination in those elements $\xi$ \autoref{not:xi} which appear in said expansion, with unknown random coefficients. Expand using \eqref{eq:exp} and equate with the expansion of step 1., thus obtaining conditions on the joint moments of the unknown random coefficients up to some inhomogeneous order.
\item The most challenging step. Solve the above moment problem by realising the random coefficients as polynomials of a low number of Gaussians. Substituting in Gaussian cubature yields a cubature formula on Wiener space.
\end{enumerate}

\subsection{Expansion of the expected signature over a symmetrised spanning
set}
{We begin by stating the main result of this subsection - the expected signature of Brownian motion expressed in representation \eqref{eq:expansion}. The remainder of the subsection will derive the result through a series of lemmas.}
\begin{proposition}\label{proposition:deg_7_expansion}
Let $B$ be a Brownian motion in $\mathbb{R}^d$. Then it has expected signature given by
\begin{equation}
\begin{aligned}\label{equation:degree_7_expansion}
&\mathbb{E}[S^{(7)}_{0,1}(\circ B)] \\
={}&1 + \sum_{1\leq i\leq d}\frac{1}{2}(\xi_i, \xi_i ) +\sum_{1\leq i,j\leq d}\Big [\frac{1}{2}(\xi_i, \xi_{ijj} ) +\frac{1}{4}(\xi_{ij}, \xi_{ij} ) + \frac{1}{8}(\xi_{i},\xi_{i},\xi_{j},\xi_{j})\Big ]\\
& + \sum_{1\leq i,j,k\leq d}\Big [\frac{1}{12}(\xi_{i},\xi_{ijjkk}) + \frac{1}{24}(\xi_{j},\xi_{iijkk}) + \frac{1}{6}(\xi_{ij},\xi_{ijkk}) + \frac{1}{12}(\xi_{ik},\xi_{ijjk})\\
& + \frac{1}{8}(\xi_{ijj},\xi_{ikk}) + \frac{1}{12}(\xi_{ijk},\xi_{ijk}) +\frac{1}{4}(\xi_{i},\xi_{i},\xi_{j},\xi_{jkk}) + \frac{1}{8}(\xi_{i},\xi_{i},\xi_{jk},\xi_{jk})\\
&+ \frac{1}{6}(\xi_{i},\xi_{j},\xi_{ik},\xi_{jk})+ \frac{1}{48}(\xi_{i},\xi_{i},\xi_{j},\xi_{j},\xi_{k},\xi_{k})\Big ].
\end{aligned}
\end{equation}
Moreover, the time-augmented case can be reduced to the above by
\begin{equation}
\begin{aligned}
\mathbb{E}[S_{0,1}^{(7)}&(\circ\widehat{B})] = \mathbb{E}[S_{0,1}^{(7)}(\circ B)] + \xi_0 + \frac 12 (\xi_{0},\xi_0) + \frac 16(\xi_{0},\xi_0, \xi_0) + \sum_{1\leq i\leq d}\Big (\frac 12 \xi_{0ii} + \frac 12 (\xi_{0},\xi_{i}, \xi_{i})\\
& +\frac 12 (\xi_0, \xi_{0ii}) + \frac 16(\xi_{0i}, \xi_{0i}) + \frac 14 (\xi_0,\xi_0,\xi_i,\xi_i)\Big ) + \sum_{1\leq i,j\leq d} \Big (\frac{1}{12}\xi_{0iijj} + \frac{1}{24}\xi_{ii0jj}\\
&+ \frac 12(\xi_0,\xi_i,\xi_{ijj}) + \frac 14(\xi_j,\xi_j,\xi_{0ii}) + \frac 14(\xi_0,\xi_{ij},\xi_{ij}) + \frac 13 (\xi_i,\xi_{0j},\xi_{ij}) + \frac 18(\xi_0,\xi_i,\xi_i,\xi_j,\xi_j)\Big ).
\end{aligned}
\end{equation}

\begin{proof}
The expected signature truncated at degree 7 is,
\[
\begin{aligned}
\mathbb{E}[S_{[0,1]}^{(7)}(\circ \widehat{B})]=&\,\,\epsilon_0 + \frac{1}{2}\epsilon_0^2+\frac{1}{6}\epsilon_0^3+\frac{1}{2}\sum_{1\leq i\leq d}\epsilon_i^2+\frac{1}{2}\sum_{1\leq i\leq d}(\epsilon_0,\epsilon_i^2)+\frac{1}{4}\sum_{1\leq i\leq d}(\epsilon_0,\epsilon_0,\epsilon_i^2)\\
&+\frac{1}{8}\sum_{1\leq i,j\leq d}(\epsilon_i^2,\epsilon_j^2)+\frac{1}{8}\sum_{1\leq i,j\leq d}(\epsilon_0,\epsilon_i^2,\epsilon_j^2)+\frac{1}{48}\sum_{1\leq i,j,k\leq d}(\epsilon_i^2,\epsilon_j^2,\epsilon_k^2) .
\end{aligned}
\]
By recalling \autoref{not:xi} and the fact that for each \(i\), $\xi_i=\e(\epsilon_i)=\epsilon_i$,
we can obtain the following fairly trivial expansions.
\[
\sum_{1\leq i\leq d}\epsilon_i^2=\frac 12\sum_{1\leq i\leq d} (\xi_i,\xi_i ),\qquad \epsilon_0^2=(\xi_0,\xi_0),\qquad \text{and},\quad\epsilon_0^3=(\xi_0,\xi_0,\xi_0).
\]
\autoref{lemma:terms_in_signature_first} - \autoref{lemma:terms_in_signature_last} will subsequently expand each remaining term in the symmetrised Eulerian representation \eqref{eq:expansion}, which will conclude the proof.
\end{proof}
\end{proposition}
\begin{remark}\label{remark:symmetry}
\label{benchmark remark}\label{ref:3d}
Each symmetric product in the expansion \eqref{equation:degree_7_expansion} of the expected signature involves no more
than three distinct basis elements $\epsilon_{i}$ and the general expansion can be deduced from the three dimensional case by symmetry. We will later see that  any Lie polynomial in the support of a degree-seven cubature on Wiener space can also be written in terms of Lie monomials involving no more than three distinct basis elements $\epsilon _{i}$ and the formula follows (provided the Lie polynomials are sufficiently symmetric) by symmetry from the three-dimensional case.
\end{remark}

Before stating and proving the expansions of each term using \eqref{eq:expansion}, we derive a general result (\autoref{cor:parity}) which reduces the required computation approximately by half.

\begin{lemma}\label{lem:reversal}
The reversal operator $*$ \eqref{eq:reversal} is a coalgebra morphism w.r.t.\ $\Delta_\shuffle$.
\begin{proof}
Let $v_1,\ldots,v_n \in V$, $x = v_1\ldots v_n$, $r(i) \coloneqq n - i + 1$. Then, recalling the notation \eqref{eq:coproduct} 
\begin{align*}
    (* \otimes *) \circ \Delta_\shuffle x 
    ={}& \sum_{I \sqcup J = [n]} (v_I)^* \otimes (v_J)^* \\
    ={}& \sum_{I \sqcup J = [n]} (v^*)_{r(I)} \otimes (v^*)_{r(J)} \\
    ={}& \sum_{I' \sqcup J' = [n]} (v^*)_{I'} \otimes (v^*)_{J'} \\
    ={}&\Delta_\shuffle(v^*) . \qedhere
\end{align*}
\end{proof}
\end{lemma}
We will say that $x \in T(V)$ is a \emph{palindrome} if $x^* = x$. Every symmetric tensor---importantly, this includes the expected signature of Brownian motion---is a palindrome, although the converse is not necessarily true; in particular, the only single words that are symmetric are tensor powers of a single letter, but many other single word palindromes exist. The invariance of palindromes under the reversal operator, combined with a parity argument originating from \autoref{prop:reversal} yields the following general result.
\begin{theorem}\label{thm:odd_level_products}
Let $m_i \geq 1$ and $x_i \in V^{\otimes m_i}$ be palindromes for $i = 1,\ldots,n$. Let $k \geq 1$ with $k \not \equiv \sum_{i = 1}^n m_i \pmod 2$. Then
\[
\e^{\star k}((x_1,\ldots,x_n)) = 0 .
\]
\begin{proof}
By symmetry and the fact that each $x_i$ is a palindrome,
\begin{align*}
(x_1,\ldots,x_n) &= \frac{1}{2 \cdot n!}\sum_{\sigma \in \mathfrak S_n} \big[ x_{\sigma(1)} \ldots x_{\sigma(n)} + x_{\sigma(n)} \ldots x_{\sigma(1)} \big] \\
&= \frac{1}{2 \cdot n!}\sum_{\sigma \in \mathfrak S_n} \big[ x_{\sigma(1)} \ldots x_{\sigma(n)} + (x_{\sigma(1)} \ldots x_{\sigma(n)})^* \big] .
\end{align*}
For $y \in T(V)$ and $f \in \mathrm{End}(T(V))$, using Sweedler notation for the convolution power of $f$ \eqref{eq:star} and the fact that reversal is a coalgebra morphism \autoref{lem:reversal}, we can write
\[
f^{\star k}(y^*) = (f((y^*)_{(1)}), \ldots, f((y^*)_{(k)})) = (f((y_{(1)})^*), \ldots, f((y_{(k)})^*)) .
\]
Let $y = (x_1,\ldots,x_n)$ and $m = \sum_{i = 1}^n m_i$. Combining the above two identities and applying the reversal property \autoref{prop:reversal} concludes the proof:
\begin{align*}
    \e^{\star k}(y + y^*) &= (\e(y_{(1)}), \ldots, \e(y_{(k)})) + (\e((y_{(1)})^*), \ldots, \e((y_{(k)})^*)) \\
    &= (1 + (-1)^{m - k})\e^{\star k}(y)\\ &= 0 . \qedhere
\end{align*}
\end{proof}
\end{theorem}

As an immediate corollary, we have the following, which can be applied to any term in the expected signature for Brownian motion (augmented with time).

\begin{corollary}\label{cor:parity}
Let $k \geq 1$ with $k \not \equiv m \pmod 2$. Then for any $i_1,\ldots,i_n \in \{1,\ldots,d\}$
\[
\e^{\star k}((\underbrace{\epsilon_0,\ldots, \epsilon_0}_{m \emph{ times}}, \epsilon_{i_1}^2,\ldots,\epsilon_{i_n}^2)) = 0 .
\]
\end{corollary}

{
We are now ready to expand each term of the Brownian motion expected signature in the symmetrised representation. In each of these \autoref{lemma:terms_in_signature_first} - \autoref{lemma:terms_in_signature_last}, the proof strategy follows the same three steps:
\begin{enumerate}
\item Begin with the representation \eqref{eq:expansion} and simplify using \autoref{cor:parity}. For tensors of even (resp. odd) length, all odd (resp. even) order symmetric products will vanish.
\item For each remaining non-zero level, compute the \(n\)-fold Eulerian idempotent by applying \eqref{eq:star} (and using \autoref{rem:cocommutative} to immediately collect terms into symmetric products). Ignore any symmetric products which include the trivially zero term \(\xi_{ii}=0\) for any \(i\).
\item Eliminate linear redundancies by combining terms which are linearly dependent. This is usually achieved by applying one of the two symmetries given in \autoref{subsec:prelims}: the reversal property \autoref{prop:reversal} and the symmetric property \autoref{prop:symmetric}. Subsequently re-index such terms to restore the ordering of the indices allowing them to be combined.
\end{enumerate}
We order Lemmas \autoref{lemma:terms_in_signature_first} - \autoref{lemma:terms_in_signature_last} in increasing order of the complexity of the expansion, beginning with the (inhomogeneous) degree-4 terms. Similar expansions for these can be compared with work on degree-\(5\) cubature constructions. The first degree-4 term is \((\epsilon_0,\epsilon_i^2)\).
\begin{lemma}\label{lemma:terms_in_signature_first}
\[
\sum_{1\leq i\leq d}(\epsilon_0,\epsilon_i^2)=\sum_{1\leq i\leq d}\Big (\xi_{0ii}+(\xi_0,\xi_i,\xi_i)\Big ) .
\]
\begin{proof}
Representation \eqref{eq:expansion} combined with \autoref{cor:parity} gives,
\[
\begin{aligned}
2\sum_{1\leq i\leq d}(\epsilon_0,\epsilon_i^2)&=\e^{\star 1}\Big (2\sum_{1\leq i\leq d}(\epsilon_0,\epsilon_i^2)\Big )+\e^{\star 3}\Big (2\sum_{1\leq i\leq d}(\epsilon_0,\epsilon_i^2)\Big )\\
&=\e^{\star 1}\Big (\sum_{1\leq i\leq d}(\epsilon_0\epsilon_i^2+\epsilon_i^2\epsilon_0)\Big )+\e^{\star 3}\Big (\sum_{1\leq i\leq d}(\epsilon_0\epsilon_i^2+\epsilon_i^2\epsilon_0)\Big )\\
&=\sum_{1\leq i\leq d}\Big (\xi_{0ii}+\xi_{ii0}\Big )+\sum_{1\leq i\leq d}2(\xi_0,\xi_i,\xi_i)\\
&=\sum_{1\leq i\leq d}\Big (2\xi_{0ii}+2(\xi_0,\xi_i,\xi_i)\Big ),
\end{aligned}
\]
as required, where we use \(\xi_{0ii}=\xi_{ii0}\) (reversal property) to obtain the final equivalence.
\end{proof}
\end{lemma}
Continuing to expand the lower order terms, the next degree-\(4\) term is $(\epsilon_i^2,\epsilon_j^2)$, which requires a slightly more involved computation.
\begin{lemma}\label{lemma:terms_in_signature_second}
\[
\sum_{1\leq i,j\leq d}(\epsilon_i^2,\epsilon_j^2)=\sum_{1\leq i,j\leq d}\Big (4(\xi_i,\xi_{ijj})+2(\xi_{ij},\xi_{ij})+(\xi_i,\xi_i,\xi_j,\xi_j)\Big ) .
\]
\begin{proof}
Representation \eqref{eq:expansion} combined with \autoref{cor:parity} gives,
\begin{equation}\label{eq:proof_ei2_ej2}
\sum_{1\leq i,j\leq d}(\epsilon_i^2,\epsilon_j^2)=\e^{\star 2}\Big (\sum_{1\leq i,j\leq d}(\epsilon_i^2,\epsilon_j^2)\Big )+\e^{\star 4}\Big (\sum_{1\leq i,j\leq d}(\epsilon_i^2,\epsilon_j^2)\Big ).
\end{equation}
Expanding the symmetric product as a pure tensor and applying linearity we get,
\[
\e^{\star k}\Big (\sum_{1\leq i,j\leq d}(\epsilon_i^2,\epsilon_j^2)\Big )=\e^{\star k}\Big (\sum_{1\leq i,j\leq d}\epsilon_i^2\epsilon_j^2\Big )=\sum_{1\leq i,j\leq d}\e^{\star k} (\epsilon_i^2\epsilon_j^2 ).
\]
For \(k=2\):
\[
\e^{\star 2}(\epsilon_i^2\epsilon_j^2)=2(\xi_i,\xi_{ijj})+2(\xi_j,\xi_{iij})+2(\xi_{ij},\xi_{ij}) .
\]
Substituting \(\xi_{iij}=\xi_{jii}\) (reversal property) before re-indexing the middle term \((i,j)\mapsto (j,i)\) we get,
\begin{equation}\label{eq:k2_iijj}
\begin{aligned}
\sum_{1\leq i,j\leq d}\e^{\star 2}(\epsilon_i^2\epsilon_j^2)&=\sum_{1\leq i,j\leq d}\Big (2(\xi_i,\xi_{ijj})+2(\xi_j,\xi_{jii})+2(\xi_{ij},\xi_{ij})\Big )\\
&=\sum_{1\leq i,j\leq d}\Big (4(\xi_i,\xi_{ijj})+2(\xi_{ij},\xi_{ij})\Big ) .
\end{aligned}
\end{equation}
For \(k=4\):
\begin{equation}\label{eq:k4_iijj}
\e^{\star 4}(\epsilon_i^2\epsilon_j^2)=(\xi_i,\xi_i,\xi_j,\xi_j) .
\end{equation}
Substituting \eqref{eq:k2_iijj} and \eqref{eq:k4_iijj} into \eqref{eq:proof_ei2_ej2} we obtain the required result.
\end{proof}
\end{lemma}
Next, we develop expansions for the degree-\(6\) terms in the expected signature. It is these terms in particular that have extremely convoluted expressions when expressed in a Hall basis, that are infeasible to compute by hand. Our representation \eqref{eq:expansion}, on the other hand, gives succinct derivations and results. First, the term \((\epsilon_0,\epsilon_0,\epsilon_i^2)\).
\begin{lemma}\label{lemma:terms_in_signature_third}
\[
\sum_{1\leq i\leq d}(\epsilon_0,\epsilon_0,\epsilon_i^2)=\sum_{1\leq i\leq d}((\xi_{0},\xi_{0ii})+\frac{2}{3}(\xi_{0i},\xi_{0i})+(\xi_{0},\xi_{0},\xi_{i},\xi_{i})\Big ) .
\]
\begin{proof}
Representation \eqref{eq:expansion} combined with \autoref{cor:parity} gives,
\begin{equation}\label{eq:0_0_ii}
\begin{aligned}
3\sum_{1\leq i\leq d}(\epsilon_0,\epsilon_0,\epsilon_i^2)&=\e^{\star 2}\Big (3\sum_{1\leq i\leq d}(\epsilon_0,\epsilon_0,\epsilon_i^2)\Big )+\e^{\star 4}\Big (3\sum_{1\leq i\leq d}(\epsilon_0,\epsilon_0,\epsilon_i^2)\Big )\\
&=\e^{\star 2}\Big (\sum_{1\leq i\leq d}(\epsilon_0^2\epsilon_i^2+\epsilon_0\epsilon_i^2\epsilon_0+\epsilon_i^2\epsilon_0^2)\Big )+\e^{\star 4}\Big (\sum_{1\leq i\leq d}(\epsilon_0^2\epsilon_i^2+\epsilon_0\epsilon_i^2\epsilon_0+\epsilon_i^2\epsilon_0^2)\Big ) .
\end{aligned}
\end{equation}
The first term simplifies by substituting \(\xi_{0i}=-\xi_{i0}\), and \(\xi_{0ii}=\xi_{ii0}\) (reversal property),
\begin{equation}\label{eq:k2_0_0_ii}
\begin{aligned}
\e^{\star 2}\Big (\sum_{1\leq i\leq d}(\epsilon_0^2\epsilon_i^2+\epsilon_0\epsilon_i^2\epsilon_0+\epsilon_i^2\epsilon_0^2)\Big )=\sum_{1\leq i\leq d}&\Big (3(\xi_{0},\xi_{0ii})+3(\xi_{0},\xi_{ii0})+2(\xi_{0i},\xi_{0i})\\
&+2(\xi_{0i},\xi_{i0})+2(\xi_{i0},\xi_{i0})\Big )\\
=\sum_{1\leq i\leq d}&\Big (6(\xi_{0},\xi_{0ii})+2(\xi_{0i},\xi_{0i})\Big ) .
\end{aligned}
\end{equation}
The second term is evaluated as,
\begin{equation}\label{eq:k4_0_0_ii}
\e^{\star 4}\Big (\sum_{1\leq i\leq d}(\epsilon_0^2\epsilon_i^2+\epsilon_0\epsilon_i^2\epsilon_0+\epsilon_i^2\epsilon_0^2)\Big )=\sum_{1\leq i\leq d}3(\xi_{0},\xi_{0},\xi_{i},\xi_{i}).
\end{equation}
Substituting \eqref{eq:k2_0_0_ii} and \eqref{eq:k4_0_0_ii} into \eqref{eq:0_0_ii} we obtain the required result.
\end{proof}
\end{lemma}
Continuing with the degree-\(6\) terms, the final two lemmas while conceptually straightforward each require a somewhat tedious calculation. Firstly, the term \((\epsilon_0,\epsilon_i^2,\epsilon_j^2)\), for which we defer the proof to \autoref{app:lemma_proofs}.
\begin{lemma}\label{lem:technical1}
\[
\begin{aligned}
\sum_{1\leq i,j\leq d}(\epsilon_0,\epsilon_i^2,\epsilon_j^2)=\sum_{1\leq i,j\leq d}&\Big (\frac{2}{3}\xi_{0iijj}+\frac{1}{3}\xi_{ii0jj}+2(\xi_0,\xi_i,\xi_{ijj})+2(\xi_i,\xi_i,\xi_{0jj})\\
&+2(\xi_{0},\xi_{ij},\xi_{ij})+ \frac{8}{3}(\xi_i,\xi_{0j},\xi_{ij})+(\xi_{0},\xi_i,\xi_i,\xi_j,\xi_j)\Big ).
\end{aligned}
\]
\end{lemma}
Finally, the expansion of \((\epsilon_i^2,\epsilon_j^2,\epsilon_k^2)\), for which again we defer the proof to \autoref{app:lemma_proofs}. This expansion is the most involved as it concerns all of the terms in three distinct basis variables.
\begin{lemma}\label{lemma:terms_in_signature_last}
\[
\begin{aligned}
\sum_{1\leq i,j,k\leq d}(\epsilon_i^2,\epsilon_j^2,\epsilon_k^2)=\sum_{1\leq i,j,k\leq d}\Big (4&(\xi_i,\xi_{ijjkk})+2(\xi_j,\xi_{iijkk})+8(\xi_{ij},\xi_{ijkk})+4(\xi_{ik},\xi_{ijjk})&
\\        &+6(\xi_{ijj},\xi_{ikk})+4(\xi_{ijk},\xi_{ijk})+12(\xi_i,\xi_k,\xi_k,\xi_{ijj})\\
&+6(\xi_i,\xi_i,\xi_{jk},\xi_{jk})+8(\xi_{i},\xi_{j},\xi_{ik},\xi_{jk})+(\xi_i,\xi_i,\xi_j,\xi_j,\xi_k,\xi_k)\Big ).
\end{aligned}
\]
\end{lemma}
}
\subsection{Degree-seven cubature on Wiener space for arbitrary dimensions}
Aided by our symmetrised representation of the degree seven truncated expected signature of Brownian motion \autoref{proposition:deg_7_expansion}, we can now construct a cubature formula for (time-augmented) Wiener space, which constitutes the main result of this paper. The measure introduced in this section is based on an ansatz in which the Lie polynomials in its support are expressed as linear combinations of Eulerian idempotents, with random coefficients given by linear combinations of products of independent Gaussian random variables. These random variables are then realised through independent Gaussian cubature formulas of suitable degree within our cubature measure. For a discrete measure $\rho = \sum_{n=1}^N \rho_n \delta_{y^{(n)}}$ with $y^{(n)}=\left( y^{(n)}_1 , \ldots , y^{(n)}_e \right)$ we will sometimes write $(y_i,\rho )$ to highlight both the particles in the support of the measure which we will also sometimes interpret as random variables.  Our main theorem is the following.

\begin{theorem}[Degree-7 cubature formula on \(d\)-dimensional Wiener space with drift]\label{theorem:deg_7_cubature_formula}
Let $(z_i,\lambda)$, $(z_{ij},\mu)$ and $(z,\eta)$ be independent Gaussian cubature formulae of $(\text{deg}=7, \text{dim}=d)$, $(\text{deg}=3, \text{dim}=d^2)$ and $(\text{deg}=2, \text{dim}=1)$ respectively. Define,
\begin{align*}
\mathcal{L}^{(n,m,r)} \coloneqq{}& \epsilon_0 + \sum_{i}\Big (z_i^{(n)}\epsilon_i+\frac{1}{\sqrt 3}z_{ii}^{(m)}\e(\epsilon_0\epsilon_i)+\frac{1}{2}\e(\epsilon_0\epsilon_i\epsilon_i)\Big ) +\sum_{i,j}\Big[ \Big( \frac{1}{\sqrt{3}}z_i^{(n)} z_{jj}^{(m)}+ \frac{1}{\sqrt{6}}z_{ij}^{(m)}\Big) \e(\epsilon_i\epsilon_j)\\
&+\frac{1}{2} z_i^{(n)}\e(\epsilon_i\epsilon_j\epsilon_j)+\frac{1}{12}\e(\epsilon_0\epsilon_i\epsilon_i\epsilon_j\epsilon_j)+\frac{1}{24}\e(\epsilon_i\epsilon_i\epsilon_0\epsilon_j\epsilon_j)\Big] +\sum_{i,j,k} \Big (\frac{1}{\sqrt{6}} z_{ij}^{(m)}z_k^{(n)}z^{(r)}\e(\epsilon_i\epsilon_j\epsilon_k)\\
&+\frac{1}{2\sqrt{3}}z_{i}^{(n)}z_{jj}^{(m)}\e(\epsilon_i\epsilon_j\epsilon_k\epsilon_k)+\frac{1}{4\sqrt{3}}z_{i}^{(n)}z_{kk}^{(m)}\e(\epsilon_i\epsilon_j\epsilon_j\epsilon_k)+\frac{1}{12}z_i^{(n)} \e(\epsilon_i\epsilon_j\epsilon_j\epsilon_k\epsilon_k)\\
&+\frac{1}{24}z_j^{(n)} \e(\epsilon_i\epsilon_i\epsilon_j\epsilon_k\epsilon_k)\Big ).
\end{align*}
and $\theta_{n,m,r} = \lambda_n\mu_m\eta_r$ for each $(n,m,r)$ that indexes the product measure $(z_i,\lambda)\times (z_{ij},\mu)\times (z,\eta)$. Then
$\mathcal{L}^{(n,m,r)}$ and $\theta_{n,m,r}$ define the Lie polynomials and weights respectively of a degree-\(7\) cubature formula on \(d\)-dimensional Wiener space with drift.

\end{theorem}
Before turning to the proof of \autoref{theorem:deg_7_cubature_formula}, we take a moment to explain how such a cubature formula is constructed. \autoref{proposition:deg_7_expansion} provides a sparse representation of the expected signature of Brownian motion in terms of symmetric products of Eulerian idempotents  (including the basis vectors $\epsilon_i$). Our aim is to first construct a measure supported on Lie polynomials, expressed as linear combinations of these Eulerian idempotents with random coefficients that match the cubature on Wiener space property. More precisely, the cubature condition \autoref{cubfreeliedef} requires that the expectation $\mathbb{E}_{\mathcal{L}_{n,m,r}} \left( \exp(L) \right)$ coincide with the expected signature of Brownian motion up to degree seven (in the inhomogeneous grading which counts drift vectors $\epsilon_0$ twice). Expanding the exponential and applying \eqref{eq:exp} reduces this to a moment-matching problem: the coefficients of the symmetric products of Eulerian idempotents, which are powers of the random Lie polynomial coefficients, must match those arising in the Brownian case.

To make the problem tractable, we adopt an ansatz in which the random coefficients are products of factors with (marginal) Gaussian distributions. This choice (i) preserves the symmetry of the Brownian expected signature, ensuring that terms which vanish in the Brownian case also vanish in our construction, and (ii) reduces the remaining moment-matching constraints to a manageable number. Additionally, the coefficients of the basis vectors $\epsilon_i$ always match the moments of a standard Gaussian up to degree seven, as seen by mapping the cubature property into the commutative algebra (cf. \cite{LV04}{, Proposition 5.1}).

After making this choice, there remains a (manageable) number of non-zero moment constraints, which lead to systems of polynomial equations for the higher-order Eulerian idempotent coefficients. To obtain real solutions, we extend the coefficient ansatz with linear combinations of products of Gaussian variables, thereby introducing enough degrees of freedom to solve these systems (compare the Lie polynomials in \autoref{theorem:deg_7_cubature_formula}).

Since the cubature on Wiener space property depends in our formulation only on the moments of the random variables, Gaussian variables may be replaced by cubature measures of appropriate degree. By taking the auxiliary Gaussian random variables to be independent, each is required to match moments only up to degree three. To see this observe that the auxiliary variables are associated to homogeneous Lie polynomials of degree at least two. At most three can occur in any symmetric product in the expansion of the cubature property up to degree seven. This observation allows the use of degree-three Gaussian cubatures for the auxiliary variables, greatly reducing the support size of the resulting cubature measure. Independent cubatures are obtained by product constructions.

The following remark illustrates how this approach is natural and simplifies the problem beyond the symmetries that force certain terms to vanish (compare also Litterer \cite{litterer} and Shinozaki \cite{shinozaki}). It explores the construction of the random coefficient for \(\e(\epsilon_i\epsilon_j)\), which is involved in several non-linear constraints. In the following $\mathbb{E}$ will denote the expectation with respect to the discrete product measure $(z_i,\lambda)\times (z_{ij},\mu)\times (z,\eta).$

\begin{remark}
Our ansatz in \autoref{theorem:deg_7_cubature_formula} suggests the following form for the coefficient of \(\e(\epsilon_i\epsilon_j)\),
\[
c_{ij} z_i z_{jj} + \hat{c}_{ij} z_{ij}.
\]
Here $c_{ij}$ and $\hat{c}_{ij}$ are unknown constants to be determined by all constraints arising from matching coefficients of the symmetric products involving $\e(\epsilon_i\epsilon_j)$ in the expansion of the expectation of the Brownian signature. From the expansion of the expectation of the Brownian signature \autoref{proposition:deg_7_expansion}, expanding the exponential in $\mathbb{E}_{\mathcal{L}_{n,m,r}} \left( \exp(L) \right)$ and equating the coefficients of the symmetric product $\big(\e(\epsilon_i\epsilon_j),\e(\epsilon_i\epsilon_j)\big)$ yields the following constraint,
\begin{equation}
\frac{1}{4}
= \frac{1}{2!}\,\mathbb{E}\!\left(\left(c_{ij}^{\,} z_i z_{jj} + \hat{c}_{ij} z_{ij}\right)^2\right)
= \frac{1}{2}\left[ c_{ij}^2\,\mathbb{E}\!\left(z_i^2 z_{jj}^2\right) + \hat{c}_{ij}^2\,\mathbb{E}\!\left(z_{ij}^2\right) \right]
= \frac{c_{ij}^2 + \hat{c}_{ij}^2}{2}.
\end{equation}

Similar equations arise from matching the coefficients of the terms $\left(\xi_i,\xi_j,\xi_{ik},\xi_{jk}\right)$, $\left(\xi_0,\xi_{ij},\xi_{ij}\right)$, $\left(\xi_i,\xi_{0j},\xi_{ij}\right)$, and $\left(\xi_i,\xi_i,\xi_{jk},\xi_{jk}\right)$ in the expansion of the expected signature,  yielding a system of polynomial equations that appears overdetermined.

However, due to the choice of ansatz several constraints are equivalent. For example, the coefficient of $\left(\xi_i,\xi_i,\xi_{jk},\xi_{jk}\right)$. By \eqref{eq:exp} and \autoref{proposition:deg_7_expansion} and recalling that for any Gaussian cubature measure we have $\mathbb{E}[z_i^2] = 1$, then,
\begin{equation}
\frac{1}{8}
= \frac{1}{2!2!}\,\mathbb{E}\!\left( z_i^2 \left(c_{jk}^{\,} z_j z_{kk} + \hat{c}_{jk} z_{jk}\right)^2 \right)
= \frac{c_{jk}^2 + \hat{c}_{jk}^2}{4}.
\end{equation}

Under the ansatz (which can be solved in the three-dimensional case and then generalised by symmetry) the complete system of equations has a real solution $c_{ij} = \frac{1}{\sqrt{3}}$ and $\hat{c}_{ij} = \frac{1}{\sqrt{6}}$. Without this internal consistency arising from the interplay of the symmetries of the expectation of the Brownian signature and the structure of our ansatz, a construction that works for arbitrary degree-seven cubatures (rather than one specific measure for the coefficients of the $\epsilon_i$) and any independent Gaussian auxiliary cubatures would not be possible.
\end{remark}

We will now present proof of our cubature formula on Wiener space by systematically verifying all moment conditions for the coefficients of the Lie polynomials in the support of our cubature formula.

{
\begin{proof}[Proof of \autoref{theorem:deg_7_cubature_formula}]
The cubature property \autoref{cubfreeliedef} can be verified by direct computation: exponentiate the Lie polynomials in the support of the cubature measure, expand in terms of symmetric products of Eulerian idempotents using \eqref{eq:exp} and compare the resulting coefficients with those in \autoref{proposition:deg_7_expansion}. Because the drift $\epsilon_0$ has coefficient one, the constraints arising from symmetric products involving only the basis vectors $\epsilon_0,\ldots,\epsilon_d$ reduce to Gaussian moment identities after mapping the cubature-on-Wiener-space property to the symmetric (commutative) algebra; equivalently, they are satisfied if and only if the coefficients $z_i$, $i=1,\ldots,d$ realise a degree-seven cubature formula (cf.~\cite{LV04}{Proposition 5.1}). Hence, the coefficients of the basis vectors $\epsilon_i$ always match the moments of a standard Gaussian up to degree seven. It therefore remains only to verify the coefficients of those symmetric products that contain at least one higher-order Lie polynomial. To keep the calculations concise, we organise this verification into three parts: terms that vanish, terms that are non-vanishing but particularly involved (these typically include lower-order Eulerian idempotents, which, since we truncate at degree seven, appear in multiple conditions), and finally all remaining non-vanishing terms that, whilst tedious, are a relatively straightforward bookkeeping exercise.

Firstly, the zero terms. Any symmetric product that does not contain an even number of instances of each basis variable $\epsilon_i$ should be zero. This holds since in the proposed cubature formula any instance of a single basis variable is matched with a corresponding Gaussian coefficient. To give two demonstrative examples, the coefficient of \(\xi_{ijj}\) features a single \(i\) which is matched by a coefficient proportional to \(z_i^n\). The coefficient of \(\xi_{ij}\) features both a single \(i\) and \(j\). These are matched by two separate coefficients, one being \(z_i^nz_{jj}^m\) and the second being \(z_{ij}^m\), both of which match a single \(i\) and single \(j\). Because the Lie polynomials feature this structure, then after exponentiation any term which contains an ``odd'' number of instances of any particular basis variable will have at least one ``odd'' Gaussian coefficient - which in expectation is always zero as required. There are also a handful of additional terms which despite containing only even instances of basis variables, do not appear in the expansion of \autoref{proposition:deg_7_expansion}. Listing these exhaustively, we can verify that the auxiliary Gaussian cubature coefficients have been carefully selected to correctly ensure these are also removed via an odd-degree Gaussian coefficient after exponentiation. The details are summarised in the following table, where all constants have been dropped for brevity.
\begin{center}
\begin{tabular}{|c|c|}
\hline
&\\[-10pt]
\textbf{Basis Term} & \textbf{Exp. Coeff.}\\
\hline
&\\[-8pt]
\((\xi_i,\xi_{0i})\) & \(\mathbb{E}[z_iz_{ii}]\)\\[4pt]
\((\xi_{0i},\xi_{ijj})\) & \(\mathbb{E}[z_{ii}z_i]\)\\[4pt]
\((\xi_i,\xi_j,\xi_{ij})\) & \(\mathbb{E}[z_iz_jz_{ij}]\)\\[4pt]
\((\xi_i,\xi_j,\xi_{ijkk})\) & \(\mathbb{E}[z_i^2z_jz_{jj}]\)\\[4pt]
\((\xi_i,\xi_k,\xi_{ijjk})\) & \(\mathbb{E}[z_i^2z_kz_{kk}]\)\\[4pt]
\((\xi_i,\xi_{jk},\xi_{ijk})\) & \(\propto \mathbb{E}[z]\)\footnotemark\\[4pt]
\hline
\end{tabular}\footnotetext{The full coefficient here is convoluted and it is simpler to observe that \(z\) only ever appears in the \(\xi_{ijk}\) coefficient.}\,\,\,
\begin{tabular}{|c|c|}
\hline
&\\[-10pt]
\textbf{Basis Term} & \textbf{Exponentiated Coefficient}\\
\hline
&\\[-8pt]
\((\xi_i,\xi_{ij},\xi_{jkk})\) & \(\mathbb{E}[z_i(z_iz_{jj}+z_{ij})(z_j+z_{jk}z_kz)]\)\\[4pt]
\((\xi_{ij},\xi_{jk},\xi_{ik})\) & \(\mathbb{E}[(z_iz_{jj}+z_{ij})(z_{j}z_{kk}+z_{jk})(z_{i}z_{kk}+z_{jk})]\)\\[4pt]
\((\xi_{0},\xi_{i},\xi_{j},\xi_{ij})\) & \(\mathbb{E}[z_iz_j(z_iz_{jj}+z_{ij})]\)\\[4pt]
\((\xi_i,\xi_j,\xi_j,\xi_{0i})\) & \(\mathbb{E}[z_iz_j^2z_{ii}]\)\\[4pt]
\((\xi_i,\xi_j,\xi_k,\xi_{ijk})\) & \(\mathbb{E}[z_iz_jz_k^2z_{ij}z]\)\\[4pt]
\((\xi_i,\xi_j,\xi_k,\xi_k,\xi_{ij})\) & \(\mathbb{E}[z_iz_jz_k^2z_{ij}]\)\\[4pt]
\hline
\end{tabular}
\end{center}

Next, the non-zero terms. As we know, the expected signature represented in \autoref{proposition:deg_7_expansion} contains a number of linear redundancies, for example \(\e(\epsilon_1\epsilon_2)=-\e(\epsilon_2\epsilon_1)\), yet both terms appear in the expansion. Our cubature formula features the same redundancies to preserve the majority of the symmetries. In this regard, for each term we aim to show that the cubature produces the correct coefficient for that exact term. There are, however, two special cases for which this process does not work and a more detailed computation must be carried out and so we present these first.
\begin{enumerate}
\item \((\xi_{ijj},\xi_{ijj})\). This term requires careful consider since it has contributions from \((\xi_{ijj},\xi_{ikk})\) and \((\xi_{ijk},\xi_{ijk})\) (both whenever \(j=k\)). The required coefficient from \autoref{proposition:deg_7_expansion} is,
\[
\frac{1}{8} + \frac{1}{12}=\frac{5}{24}.
\]
The total coefficient after exponentiating the cubature polynomials is,
\[
\frac{1}{2!}\cdot \Big (\frac{1}{4}\mathbb{E}[z_i^2] + \frac{1}{6}\mathbb{E}[z_{ij}^2z_{j}^2z^2]\Big )=\frac{5}{24},
\]
as required, where \(\frac{1}{2!}\) is the coefficient from the second level of the \(\exp\) function, and the \(2\) terms inside the brackets are the two ways to obtain \((\xi_{ijj},\xi_{ijj})\) by combining the terms of the cubature Lie polynomial.
\item \((\xi_i,\xi_i,\xi_{ij},\xi_{ij})\) (where \(i <j\)). This term is particularly involved due to the asymmetry of the \(\xi_{ij}\) coefficient. The required coefficient from \autoref{proposition:deg_7_expansion} (noting that there are positive contributions from both \((\xi_i,\xi_i,\xi_{ij},\xi_{ij})\) and \((\xi_i,\xi_i,\xi_{ji},\xi_{ji})\) as well as \((\xi_i,\xi_k,\xi_{ij},\xi_{kj})\) when \(k=i\)) is,
\[
\frac{1}{8} + \frac{1}{8} + \frac{1}{6} = \frac{5}{12}.
\]
The total coefficient after exponentiating the cubature polynomials is,
\[
\binom{4}{2}\cdot \frac{1}{4!}\cdot \Big (\frac{1}{3}\mathbb{E}[z_i^4z_{jj}^2] + \frac{1}{3}\mathbb{E}[z_i^2z_{ii}^2z_j^2]+\frac{1}{6}\mathbb{E}[z_i^2z_{ij}^2]+\frac{1}{6}\mathbb{E}[z_i^2z_{ji}^2]\Big )=\frac{5}{12},
\]
as required, where \(\binom{4}{2}\) is the size of the permutation group of \(\{a,a,b,b\}\), \(\frac{1}{4!}\) is the coefficient from the fourth level of the \(\exp\) function, and the \(4\) terms inside the brackets are the only possible ways to obtain either \((\xi_i,\xi_i,\xi_{ij},\xi_{ij})\) or \((\xi_i,\xi_i,\xi_{ji},\xi_{ji})\) by combining the terms of the cubature Lie polynomial.
\end{enumerate}

The remaining terms require no groupings. This is not because the linear redundancies do not exist, rather the symmetries do not make them a problem and leaving redundant terms ungrouped makes for a simpler and more efficient proof. Brief details on each case are summarised in the following table, where for each non-zero term in the expansion of the the expected signature \autoref{proposition:deg_7_expansion}, 
the column ``cubature coefficient'' gives the coefficient which results from exponentiating and summing the proposed cubature formula. In each case, the coefficient is computed by applying the tensor exponential identity \eqref{eq:exp} and scaling by the (expectation of) the coefficients of each constituent term in the symmetrised product as they appear in our Wiener space cubature formula. The result can be verified as equal to the target value obtained directly from \autoref{proposition:deg_7_expansion}. Unless otherwise indicated, each row holds for any (non-trivial\footnote{By non-trivial choices, we are referring to those which make the symmetrised product term a non-zero tensor.}) choice of \(i,j,k\). If a term is marked with an asterisk this is to indicate \(i\neq j\) and if a term is marked with a dagger this is to indicate \(i,j,k\) are all distinct.
\begin{center}
\begin{tabular}{|c|c|c|}
\hline
&&\\[-10pt]
\textbf{Term} & \textbf{Cub. coeff.} & \textbf{\ref{proposition:deg_7_expansion}} \\
\hline
&&\\[-8pt]
\(\xi_{0ii}\) & \(\frac{1}{2}\) & \(\frac{1}{2}\) \\[4pt]
\(\xi_{0iijj}\) & \(\frac{1}{12}\) & \(\frac{1}{12}\) \\[4pt]
\(\xi_{ii0jj}\) & \(\frac{1}{24}\) & \(\frac{1}{24}\) \\[4pt]
\((\xi_0,\xi_{0ii})\) & \(\frac{2}{2!}(\frac{1}{2})\) & \(\frac{1}{2}\) \\[4pt]
\((\xi_{0i},\xi_{0i})\) & \(\frac{1}{2!} \mathbb{E}[\frac{1}{3}z_{ii}^2]\) & \(\frac{1}{6}\)\\[4pt]
\((\xi_i,\xi_{ijj})\) & \(\frac{2}{2!} \mathbb{E}[\frac{1}{2}z_i^2] \) & \(\frac 12\) \\[4pt]
\((\xi_{i},\xi_{ijjkk})\) & \(\frac{2}{2!} \mathbb{E}[\frac{1}{12}z_{i}^2]\) & \(\frac{1}{12}\) \\[4pt]
\((\xi_{j},\xi_{iijkk})\) & \(\frac{2}{2!} \mathbb{E}[\frac{1}{24}z_{i}^2]\) & \(\frac{1}{24}\) \\[4pt]
\((\xi_{ij},\xi_{ijkk})\) & \(\frac{2}{2!}\mathbb{E}[\frac{1}{6}z_{i}^2z_{jj}^2]\) & \(\frac{1}{6}\) \\[4pt]
\((\xi_{ik},\xi_{ijjk})\) & \(\frac{2}{2!}\mathbb{E}[\frac{1}{12}z_{i}^2z_{kk}^2]\) & \(\frac{1}{12}\) \\[4pt]
\((\xi_{ijk},\xi_{ijk})\) & \(\frac{1}{2!}\mathbb{E}[\frac{1}{6}z_{ij}^2z_k^2z^2]\) & \(\frac{1}{12}\) \\[4pt]
\hline
\end{tabular}\,\,\,\begin{tabular}{|c|c|c|}
\hline
&&\\[-10pt]
\textbf{Term} & \textbf{Cubature coefficient} & \textbf{\ref{proposition:deg_7_expansion}} \\
\hline
&&\\[-8pt]
\((\xi_{ijj},\xi_{ikk})^\dagger\) & \(\frac{1}{2!}\mathbb{E}[\frac{1}{4}z_{i}^2]\) & \(\frac{1}{8}\) \\[4pt]\((\xi_{ij},\xi_{ij})\) & \(\frac{1}{2!} \mathbb{E}[\frac{1}{3}z_i^2z_{jj}^2 +\frac 16 z_{ij}^2]\) & \(\frac 14\) \\[4pt]
\((\xi_0,\xi_i,\xi_{ijj})\) & \(\frac{6}{3!}\mathbb{E}[\frac{1}{2} z_i^2]\) & \(\frac{1}{2}\)\\[4pt]
\((\xi_i,\xi_i,\xi_{0jj})\) & \(\frac{3}{3!}\mathbb{E}[\frac{1}{2}z_i^2]\) & \(\frac{1}{4}\)\\[4pt]
\((\xi_0,\xi_{ij},\xi_{ij})\) & \(\frac{3}{3!}\mathbb{E}[\frac{1}{3}z_i^2z_{jj}^2+\frac{1}{6}z_{ij}^2]\) & \(\frac{1}{4}\)\\[4pt]
\((\xi_i,\xi_{0j},\xi_{ij})\) & \(\frac{6}{3!}\mathbb{E}[\frac{1}{3}z_i^2z_{jj}^2]\) & \(\frac{1}{3}\)\\[4pt]
\((\xi_{i},\xi_i,\xi_i,\xi_{ikk})\) & \(\frac{4}{4!}\mathbb{E}[\frac{1}{2}z_{i}^4]\) & \(\frac{1}{4}\) \\[4pt]
\((\xi_{i},\xi_i,\xi_j,\xi_{jkk})^*\) & \(\frac{12}{4!}\mathbb{E}[\frac{1}{2}z_{i}^2z_j^2]\) & \(\frac{1}{4}\) \\[4pt]
\((\xi_{i},\xi_i,\xi_{jk},\xi_{jk})^*\) & \(\frac{6}{4!}\mathbb{E}[\frac{1}{3}z_i^2z_j^2z_{kk}^2 +\frac 16 z_i^2z_{jk}^2]\) & \(\frac{1}{8}\) \\[4pt]
\((\xi_i,\xi_{j},\xi_{ik},\xi_{jk})^*\)& \(\frac{12}{4!}\mathbb{E}[\frac{1}{3}z_i^2z_j^2z_{kk}^2]\) & \(\frac{1}{3}\)\\[4pt]
&&\\[4pt]
\hline
\end{tabular}
\end{center}
\end{proof}
}
\begin{remark}
All proofs in this paper are self-contained. In addition, the  cubature measure constructed in \autoref{theorem:deg_7_cubature_formula} have also been verified by direct computation in Python (for dimension-\(3\), the general case follows by symmetry as discussed in \autoref{remark:symmetry}); the corresponding code is available at \cite{code}. The verification code makes use of the computational library for the free Lie algebra developed by Reizenstein \cite{FLApy}. We provide this code for readers who may wish to construct further cubature measures on path space and who may find it a useful resource.
\end{remark}
\subsection{Breaking symmetry to reduce the support size of the cubature measure further}\label{sec:brokensym}
Given that the Eulerian idempotent method produces a cubature formula with a number of redundant terms, it is natural to question whether these can be removed to reduce the size of the cubature formula. For the degree-7 cubature formula given in \autoref{theorem:deg_7_cubature_formula}, the only redundancy that directly affects the number of cubature points is the inclusion of both the terms $\e(ij)$ and $\e(ji)$ for any $i,j$, since these require $z_{ij}$ to have dimension $d^2$. These are linked by the anti-symmetry $\e(ij)=-\e(ji)$, so we could hope to (asymptotically) halve the dimension of $z_{ij}$ by only considering $i\leq j$. Indeed this can be achieved by the following corollary.
\begin{corollary}[Degree-7 cubature formula with broken symmetries] \label{cor:sym_rem}
Let $(z_i,\lambda)$, $(z_{ij},\mu)$ and $(z,\eta)$ be independent Gaussian cubature formulae of $(\text{deg}=7, \text{dim}=d)$, $(\text{deg}=3, \text{dim}=\frac{1}{2}d (d+1)$) and $(\text{deg}=2, \text{dim}=1)$ respectively, where we understand $z_{ij}$ has dimensions only for $i\leq j$. Define,
\begin{align*}
\mathcal{L}^{(n,m,r)} \coloneqq{}&\epsilon_0 + \sum_{1\leq i\leq d}\Big (z_i^{(n)}\epsilon_i+\frac{1}{\sqrt 3}z_{ii}^{(m)}\e(\epsilon_0\epsilon_i)+\frac{1}{2}\e(\epsilon_0\epsilon_i\epsilon_i)\Big ) +\sum_{1\leq i< j\leq d}\Big[ \Big( \frac{1}{\sqrt{3}}z_i^{(n)} z_{jj}^{(m)}\\
&-\frac{1}{\sqrt{3}}z_j^{(n)} z_{ii}^{(m)}+ \frac{1}{\sqrt{3}}z_{ij}^{(m)}\Big) \e(\epsilon_i\epsilon_j) +\frac{1}{2} z_i^{(n)}\e(\epsilon_i\epsilon_j\epsilon_j)+\frac{1}{2} z_j^{(n)}\e(\epsilon_j\epsilon_i\epsilon_i)\Big]\\
&+\sum_{1\leq i,j\leq d}\Big (\frac{1}{\sqrt{6}}z_{ii}^{(m)}z_j^{(n)}z^{(r)}\e(\epsilon_i\epsilon_j\epsilon_i) + \frac{1}{12}\e(\epsilon_0\epsilon_i\epsilon_i\epsilon_j\epsilon_j)+\frac{1}{24}\e(\epsilon_i\epsilon_i\epsilon_0\epsilon_j\epsilon_j)\Big )\\
&+\sum_{\substack{1\leq i<j\leq d\\1\leq k\leq d}} \Big (\frac{1}{\sqrt{3}} z_{ij}^{(m)}z_k^{(n)}z^{(r)}\e(\epsilon_i\epsilon_j\epsilon_k)+\frac{1}{2\sqrt{3}}z_i^{(n)}z_{jj}^{(m)}\big [\e(\epsilon_i\epsilon_j\epsilon_k\epsilon_k) - \e(\epsilon_j\epsilon_i\epsilon_k\epsilon_k)\\
&+\e(\epsilon_i\epsilon_k\epsilon_k\epsilon_j)\big ]\Big )+\sum_{1\leq i,j,k\leq d}\Big (\frac{1}{12}z_i^{(n)}e(\epsilon_i\epsilon_j\epsilon_j\epsilon_k\epsilon_k)+\frac{1}{24}z_j^{(n)}e(\epsilon_i\epsilon_i\epsilon_j\epsilon_k\epsilon_k)\Big ).
\end{align*}
and $\theta_{n,m,r} = \lambda_n\mu_m\eta_r$ for each $(n,m,r)$ that indexes $(z_i,\lambda)\times (z_{ij},\mu)\times (z,\xi)$. Then
$\mathcal{L}^{(n,m,r)}$ and $\theta_{n,m,r}$ define the Lie polynomials and weights respectively of a degree-\(7\) cubature formula on \(d\)-dimensional Wiener space with drift.
\begin{proof}
Almost all computations are the same as in the proof of \autoref{theorem:deg_7_cubature_formula}, with the exception of any involving \(\xi_{ij}\). The following table gives the computations for all such symmetric product terms, which again can be directly compared with the required coefficient taken from \autoref{proposition:deg_7_expansion}. Throughout we have imposed \(i<j\) and used the fact that \(\xi_{ji}=-\xi_{ij}\) to combine terms where relevant. Special attention should be given to the symmetrised term ``\((\xi_i,\xi_{j},\xi_{ik},\xi_{jk})\)'', where one should group the case \(i<j<k\) with \(k<i<j\) and separately group the case \(i<k<j\) with \(j<k<i\). The calculation for both of these groupings is the same, but for brevity only written once in the table.
\begin{center}
\begin{tabular}{|c|c|c|}
\hline
&&\\[-10pt]
\textbf{Term} & \textbf{Cub. coeff.} & \textbf{\ref{proposition:deg_7_expansion}} \\
\hline
&&\\[-8pt]
\((\xi_{ij},\xi_{ijkk})\) & \(\frac{2}{2!}\mathbb{E}[\frac{1}{6}z_{i}^2z_{jj}^2]\) & \(\frac{1}{6}\) \\[4pt]
\((\xi_{ji},\xi_{jikk})\) & \(\frac{2}{2!}\mathbb{E}[\frac{1}{6}z_{i}^2z_{jj}^2]\) & \(\frac{1}{6}\) \\[4pt]
\((\xi_{ij},\xi_{ikkj})\) & \(\frac{2}{2!}\mathbb{E}[\frac{1}{6}z_{i}^2z_{jj}^2]\) & \(\frac{1}{6}\) \\[4pt]
\((\xi_i,\xi_{0j},\xi_{ij})\) & \(\frac{6}{3!}\mathbb{E}[\frac{1}{3}z_i^2z_{jj}^2]\) & \(\frac{1}{3}\)\\[4pt]
\((\xi_i,\xi_{0j},\xi_{ji})\) & \(\frac{6}{3!}\mathbb{E}[\frac{1}{3}z_i^2z_{jj}^2]\) & \(\frac{1}{3}\)\\[4pt]
\hline
\end{tabular}\,\,\,\begin{tabular}{|c|c|c|}
\hline
&&\\[-10pt]
\textbf{Term} & \textbf{Cubature coefficient} & \textbf{\ref{proposition:deg_7_expansion}} \\
\hline
&&\\[-8pt]
\((\xi_{ij},\xi_{ij})\) & \(\frac{1}{2!}\mathbb{E}[\frac{1}{3}z_i^2z_{jj}^2 + \frac{1}{3}z_j^2z_{ii}^2 +\frac 13 z_{ij}^2]\) & \(\frac 12\)\\[4pt]
\((\xi_0,\xi_{ij},\xi_{ij})\) & \(\frac{3}{3!}\mathbb{E}[\frac{1}{3}z_i^2z_{jj}^2+ \frac{1}{3}z_j^2z_{ii}^2 +\frac 13 z_{ij}^2]\) & \(\frac{1}{2}\)\\[4pt]
\((\xi_{i},\xi_i,\xi_{ij},\xi_{ij})\) & \(\frac{6}{4!}\mathbb{E} [\frac{1}{3}z_i^4z_{jj}^2+\frac{1}{3}z_i^2z_j^2z_{ii}^2+\frac{1}{3}z_i^2z_{ij}^2]\) & \(\frac{5}{12}\)\\[1pt]
\((\xi_k,\xi_k,\xi_{ij},\xi_{ij})\) & \(\frac{6}{4!}\mathbb{E}[\frac{z_k^2}{3}(z_i^2z_{jj}^2 + z_j^2z_{ii}^2 + z_{ij}^2)]\) & \(\frac{1}{4}\) \\[4pt]
\((\xi_i,\xi_{j},\xi_{ik},\xi_{jk})\)& \(\frac{12}{4!}\mathbb{E}[\frac{1}{3}z_i^2z_j^2z_{kk}^2]\) & \(\frac{1}{6}\)\\[4pt]
\hline
\end{tabular}
\end{center}
\end{proof}
\end{corollary}

\subsection{Comparison to existing constructions} \label{sec:gaussiancubatures}
\subsubsection{Positive Gaussian degree-seven cubature measures} \label{sec:fdGaussians}
Deterministic constructions of cubature measures on Wiener space all rely on the existence of Gaussian cubatures that are exact to the same degree and positive wights. The existence of cubature measures (both Gaussian and on Wiener space) is guaranteed by Tchakaloff's theorem (see Bayer, Teichmann \cite{BT2}). Stroud \cite{Str71} provides several examples of such measures. A degree-three Gaussian cubature in $d$ dimensions can be realised by a measure with support of $2d$ particles (see Stroud \cite{Str71}, formula $E_n^{r^2}$ 3-1, p. 315). Degree-seven measures have been constructed for $d=3$ with support of $27$ points (\cite{Str71}, formula $E_3^{r^2}$ 7-1, p. 327), $d=4$ with $49$ particles (\cite{Str71}, formula $E_4^{r^2}$ 7-1, p. 329)
and $3 \leq d \leq 8$ with $2^{d+1}+4d^2$ particles (\cite{Str71}, formula $E_3^{r^2}$ 7-2, p. 319).

For higher dimensions Gaussian cubature measures with polynomial size support can be obtained by applying recombination (see
Litterer, Lyons \cite{LL12}). A more efficient, enhanced recombination algorithms is due to Tchernychova \cite{T16} who improves the efficiency of the original recombination algorithm by a full order and applies them specifically to the construction of \say{Caratheodory} cubature measures. Her algorithm iteratively extends the dimension of a $d$-dimensional cubature measure by taking the
product with a one-dimensional Gaussian quadrature measure, followed by a recombination step with respect to the polynomials of degree at most $m$ in $d+1$ dimensions. 

Let $m=2k+1$. Note that if $G$ is random variable with all even moments up to degree $2k$ matching the standard normal distribution and $\Lambda$ is an independent Bernoulli random variable then $G$ then $G \Lambda$ has standard normal moments up to degree $2k+1$. Hence, the ``Caratheodory'' cubature can be realised by recombination with respect to even polynomials and taking the product with the Bernoulli distribution. This leads to cubature measures with support of size at most  $2 \dim(G(d,2k))+1$, where $G(d,2k)$ denotes the space of even polynomials of degree at most $2k$ in $d$ variables .

Möller \cite{M79} proves lower bounds for the size of support of Gaussian cubature measures. These bounds grow as a cubic polynomial in the dimension of the underlying space and have been approached for some dimensions for degree-five measures with positive weights (see e.g.\ Victoir \cite{victoir}). However, to the best of our knowledge, no such examples are currently known for degree-seven in higher dimensions.

\subsubsection{Efficiency compared to existing constructions} \label{sec:effic}
Shinozaki's construction \cite{shinozaki} provides a general framework for degree-seven cubature measures on Wiener space that, in principle, applies to arbitrary dimensions. The moment conditions necessary for this cubature are outlined on  pp. 907-910 \cite{shinozaki}  and are somewhat more complex than in our construction. Due to its algebraic complexity, explicit constructions and implementations have only been provided for two-dimensional Brownian motion. Discrete random variables/cubature measures are not explicitly constructed and to implement the measure for two-dimensional noise Ninomiya and Shinozaki \cite{ninomiya} use a product construction based on Gauss quadrature, yielding a seven-dimensional degree-seven Gaussian cubature measure with a support size of $5^7 = 78,125$ points. The dimension of the auxiliary Gaussian measures grows quickly and for the case of three-dimensional noise the moment conditions are based on a Gaussian cubature measure on $\mathbb{R}^{18}$.
\\ \\
Making a direct comparison in this case is challenging, as \cite{shinozaki} does not focus on optimizing the support size of their measures. While a product construction of high-dimensional degree-seven Gaussian cubature is theoretically possible, it would result in a very large support. It may be possible to reduce the support size in \cite{shinozaki} for higher-dimensional noise by incorporating some of our ideas, which leverage lower-degree Gaussian auxiliary cubatures to construct discrete random variables satisfying the moment conditions. However, given the algebraic complexity of the construction, this would likely require nontrivial modifications to the formulae, which we have not attempted to make.

The randomised construction
from Hayakawa and Tanaka \cite{hayakawa-tanaka} applied to $d$-dimensional
Brownian motion leads, when computationally tractable, to a measure of degree seven with size of support bounded above by  the dimension of the free tensor 
algebra over $\mathbb{R}\oplus \mathbb{R}^{d}$ truncated at inhomogeneous degree 7 which we denote by $D_{d}^{7}$.
Note that%
\begin{equation*}
D_{d}^{m}= \mathrm{dim}\widetilde \pi_m ( T(\mathbb R^{1+d})) = \sum_{k = 0}^m \sum_{\substack{i,j \geq 0 \\ i + 2j \leq m \\ i+j = k}}
{k \choose j} d^{k-j} =
\sum_{k=0}^{m}\sum_{j=0}^{( m-k)\wedge k }\binom{k}{j}%
d^{k-j}, 
\end{equation*}%
which gives $D_{2}^{7}=696$, $D_{3}^{7}=5,632$, $D_{4}^{7}=30,348$ , $D_{5}^{7}=121,554$ and $D_{6}^{7}=392,464$.  Hence, the
randomised construction of Hayakawa-Tanaka \cite{hayakawa-tanaka} \ which is
solved e.g.\ using linear programming quickly
becomes computationally intractable. This compares with our explicit formulae
with support of size $ 2^2 \times 2 \times 12 =96$, $2 \times 3^2 \times 2 \times 27 = 972$, $2 \times 4^2 \times 2 \times 49 = 3,136$, $2 \times  5^2 \times 2 \times \left(2^6 + 4 \times 25\right) = 16,400$ and $2 \times 6^2 \times 2 \times \left(2^7 + 4 \times 36\right) = 39168$ respectively obtained from our
construction for $2,3,4,5,6$-dimensional noise using Stroud's Gaussian cubature formulae referenced in the previous subsection. The formula obtained in \autoref{sec:brokensym} by breaking some of the symmetries of the Eulerian idempotent reduces the size of the support by $d(d+1)/(2d^2)$ (by half for the special case of two-dimensional noise). 

We have summarised the results of our comparison for $d=2,3,4,5$ in the following table.

\begin{table}[h!]
\centering
\begin{tabular}{l|c|c|c|c|c|}
\cline{2-6}
\multicolumn{1}{c|}{} & \multicolumn{5}{c|}{Dimension} \\ 
\cline{2-6}
\multicolumn{1}{c|}{} & 2 & 3 & 4 & 5 & 6 \\ 
\hline
Hall-Basis $d=2$: \cite{litterer}, $d=3$: \cite{timothy}    & 48 & 648 &  &  &  \\ \hline
Moment Similar-Families: \cite{ninomiya}   & 78,125 &  &  &  &  \\ \hline
Randomised Construction: \cite{hayakawa-tanaka}   & 696 & 5,632 & 30,348 & 121,554  & 392,464 \\ \hline
Eulerian Idempotent: Theorem \ref{theorem:deg_7_cubature_formula} & 96   & 972  & 3,136 & 16,400 & 39,168 \\ 
\hline
Eulerian Idempotent: Corollary \ref{cor:sym_rem}      & 48  & 648 & 1,960 & 9,840 & 22,848 \\ 
\hline
\end{tabular}
\caption{Comparison of size of support of different constructions for degree seven.}
\label{tab:example}
\end{table}

\appendix
\section{Proofs of the technical lemmas}\label{app:lemma_proofs}
\begin{proof} [Proof of \autoref{lem:technical1}]
Representation \eqref{eq:expansion} combined with \autoref{cor:parity} gives,
\begin{equation}\label{eq:proof_e0_ei2_ej2}
3\sum_{1\leq i,j\leq d}(\epsilon_0,\epsilon_i^2,\epsilon_j^2)=\e^{\star 1}\Big (3\sum_{1\leq i,j\leq d}(\epsilon_0,\epsilon_i^2,\epsilon_j^2)\Big )+\e^{\star 3}\Big (3\sum_{1\leq i,j\leq d}(\epsilon_0,\epsilon_i^2,\epsilon_j^2)\Big )+\e^{\star 5}\Big (3\sum_{1\leq i,j\leq d}(\epsilon_0,\epsilon_i^2,\epsilon_j^2)\Big ).
\end{equation}
Expanding the symmetric product as a tensor and applying linearity we get,
\[
\e^{\star k}\Big (3\sum_{1\leq i,j\leq d}(\epsilon_0,\epsilon_i^2,\epsilon_j^2)\Big )=\sum_{1\leq i,j\leq d}\e^{\star k} (\epsilon_0\epsilon_i^2\epsilon_j^2 +\epsilon_i^2\epsilon_0\epsilon_j^2 +\epsilon_i^2\epsilon_j^2\epsilon_0).
\]
For \(k=1\), by substituting \(\xi_{iijj0}=\xi_{0jjii}\) (reversal property) and subsequently re-indexing the third term,
\begin{equation}\label{eq:k1_0iijj}
\begin{aligned}
\sum_{1\leq i,j\leq d}\e^{\star 1} (\epsilon_0\epsilon_i^2\epsilon_j^2 +\epsilon_i^2\epsilon_0\epsilon_j^2 +\epsilon_i^2\epsilon_j^2\epsilon_0)&=\sum_{1\leq i,j\leq d}(\xi_{0iijj}+\xi_{ii0jj}+\xi_{iijj0})\\
&=\sum_{1\leq i,j\leq d}(\xi_{0iijj}+\xi_{ii0jj}+\xi_{0jjii})\\
&=\sum_{1\leq i,j\leq d}(2\xi_{0iijj}+\xi_{ii0jj}) .
\end{aligned}
\end{equation}
For \(k=3\):
\begin{equation}\label{eq:k3_0iijj_raw}
\begin{aligned}
\sum_{1\leq i,j\leq d}\e^{\star 3} (\epsilon_0\epsilon_i^2\epsilon_j^2 +\epsilon_i^2\epsilon_0\epsilon_j^2 +\epsilon_i^2\epsilon_j^2\epsilon_0)=\sum_{1\leq i,j\leq d}&\Big (3(\xi_0,\xi_i,\xi_{ijj})+3(\xi_0,\xi_j,\xi_{iij})+2(\xi_i,\xi_i,\xi_{0jj})\\
&+(\xi_i,\xi_i,\xi_{jj0})+(\xi_j,\xi_j,\xi_{0ii})+2(\xi_j,\xi_j,\xi_{ii0})\\
&+4(\xi_{i},\xi_{j},\xi_{0ij}) +4(\xi_{i},\xi_{j},\xi_{i0j}) +4(\xi_{i},\xi_{j},\xi_{ij0})\\
&+6(\xi_{0},\xi_{ij},\xi_{ij})+8(\xi_i,\xi_{0j},\xi_{ij})+4(\xi_i,\xi_{j0},\xi_{ij})\\
&+ 4(\xi_j,\xi_{0i},\xi_{ij})+8(\xi_j,\xi_{i0},\xi_{ij})\Big ) .
\end{aligned} 
\end{equation}
By substituting \(\xi_{iij}=\xi_{jii}\) (reversal property) and re-indexing of the second term,
\begin{equation}\label{eq:0_i_ijj}
\sum_{1\leq i,j\leq d}\Big (3(\xi_0,\xi_i,\xi_{ijj})+3(\xi_0,\xi_j,\xi_{iij})\Big )=\sum_{1\leq i,j\leq d}6(\xi_0,\xi_i,\xi_{ijj}) .
\end{equation}
Similarly, substituting \(\xi_{ii0}=\xi_{0ii}\) (reversal property) and re-indexing appropriate terms,
\begin{equation}\label{eq:i_i_0jj}
\sum_{1\leq i,j\leq d}\Big (2(\xi_i,\xi_i,\xi_{0jj})+(\xi_i,\xi_i,\xi_{jj0})+(\xi_j,\xi_j,\xi_{0ii})+2(\xi_j,\xi_j,\xi_{ii0})\Big )=\sum_{1\leq i,j\leq d}6(\xi_i,\xi_i,\xi_{0jj}) .
\end{equation}
By re-indexing of the middle term \((i,j)\mapsto (j,i)\) before substituting \(\xi_{ij0}+\xi_{j0i}+\xi_{0ij}=0\) (reversal + symmetric property),
\begin{equation}\label{eq:i_j_0ij}
\sum_{1\leq i,j\leq d}\Big (4(\xi_{i},\xi_{j},\xi_{0ij}) +4(\xi_{i},\xi_{j},\xi_{i0j}) +4(\xi_{i},\xi_{j},\xi_{ij0})\Big )=\sum_{1\leq i,j\leq d}4(\xi_{i},\xi_{j},\xi_{0ij}+ \xi_{j0i}+\xi_{ij0} )=0 .
\end{equation}
Finally, by re-indexing the third and fourth terms \((i,j)\mapsto (j,i)\) and substituting \(\xi_{ji}=-\xi_{ij}\) and \(\xi_{j0}=-\xi_{0j}\) (reversal / symmetric property),
\begin{equation}\label{eq:i_0j_ij}
\sum_{1\leq i,j\leq d}\Big (8(\xi_i,\xi_{0j},\xi_{ij})+4(\xi_i,\xi_{j0},\xi_{ij})+ 4(\xi_j,\xi_{0i},\xi_{ij})+8(\xi_j,\xi_{i0},\xi_{ij})\Big )=\sum_{1\leq i,j\leq d}8(\xi_i,\xi_{0j},\xi_{ij}) .
\end{equation}
Substituting \eqref{eq:0_i_ijj}, \eqref{eq:i_i_0jj}, \eqref{eq:i_j_0ij} and \eqref{eq:i_0j_ij} into \eqref{eq:k3_0iijj_raw} yields,
\begin{equation}\label{eq:k3_0iijj}
\begin{aligned}
\sum_{1\leq i,j\leq d}\e^{\star 3} (\epsilon_0\epsilon_i^2\epsilon_j^2 +\epsilon_i^2\epsilon_0\epsilon_j^2 +\epsilon_i^2\epsilon_j^2\epsilon_0)=\sum_{1\leq i,j\leq d}&\Big (6(\xi_0,\xi_i,\xi_{ijj})+6(\xi_i,\xi_i,\xi_{0jj})\\
&+6(\xi_{0},\xi_{ij},\xi_{ij})+ 8(\xi_i,\xi_{0j},\xi_{ij})\Big ) .
\end{aligned} 
\end{equation}
For \(k=5\),
\begin{equation}\label{eq:k5_0iijj}
\sum_{1\leq i,j\leq d}\e^{\star 5} (\epsilon_0\epsilon_i^2\epsilon_j^2 +\epsilon_i^2\epsilon_0\epsilon_j^2 +\epsilon_i^2\epsilon_j^2\epsilon_0)=\sum_{1\leq i,j\leq d}3(\xi_{0},\xi_i,\xi_i,\xi_j,\xi_j) .
\end{equation}
Substituting \eqref{eq:k1_0iijj}, \eqref{eq:k3_0iijj} and \eqref{eq:k5_0iijj} into \eqref{eq:proof_e0_ei2_ej2} we obtain the required result. 
\end{proof}

We conclude this appendix with the proof of the second technical lemma required for the expansion of the expected signature.

\begin{proof}[Proof of \autoref{lemma:terms_in_signature_last}]
Representation \eqref{eq:expansion} combined with \autoref{cor:parity} gives,
\begin{equation}\label{eq:proof_ei2_ej2_ek2}
\sum_{1\leq i,j,k\leq d}(\epsilon_i^2,\epsilon_j^2,\epsilon_k^2)=\e^{\star 2}\Big (\sum_{1\leq i,j\leq d}(\epsilon_i^2,\epsilon_j^2,\epsilon_k^2)\Big )+\e^{\star 4}\Big (\sum_{1\leq i,j\leq d}(\epsilon_i^2,\epsilon_j^2,\epsilon_k^2)\Big )+\e^{\star 6}\Big (\sum_{1\leq i,j\leq d}(\epsilon_i^2,\epsilon_j^2,\epsilon_k^2)\Big ).
\end{equation}
Expanding the symmetric product as a pure tensor and applying linearity we get,
\[
\e^{\star k}\Big (\sum_{1\leq i,j\leq d}(\epsilon_i^2,\epsilon_j^2,\epsilon_k^2)\Big )=\e^{\star k}\Big (\sum_{1\leq i,j\leq d}\epsilon_i^2\epsilon_j^2\epsilon_k^2\Big )=\sum_{1\leq i,j\leq d}\e^{\star k} (\epsilon_i^2\epsilon_j^2\epsilon_k^2).
\]
For \(k=2\):
\begin{equation}\label{eq:i2j2k2_k2}
\begin{aligned}
\e^{\star 2}(\epsilon_i^2\epsilon_j^2\epsilon_k^2)=2(\xi_i&,\xi_{ijjkk})+2(\xi_j,\xi_{iijkk})+2(\xi_k,\xi_{iijjk})+4(\xi_{ij},\xi_{ijkk})\\
&+4(\xi_{ik},\xi_{ijjk})+4(\xi_{jk},\xi_{iijk})+2(\xi_{iij},\xi_{jkk})\\
&+2(\xi_{ijj},\xi_{ikk})+2(\xi_{iik},\xi_{jjk})+4(\xi_{ijk},\xi_{ijk}) .
\end{aligned}
\end{equation}
Substituting \(\xi_{iijjk}=\xi_{kjjii}\) (reversal property) and re-indexing the latter term by \((i,j,k)\mapsto (k,j,i)\) we obtain,
\begin{equation}\label{eq:i_ijjkk}
\begin{aligned}
\sum_{1\leq i,j,k\leq d}\Big (2(\xi_i,\xi_{ijjkk})+2(\xi_k,\xi_{iijjk})\Big )&=\sum_{1\leq i,j,k\leq d}\Big (2(\xi_i,\xi_{ijjkk})+2(\xi_k,\xi_{kjjii})\Big )\\
&=\sum_{1\leq i,j,k\leq d}4(\xi_i,\xi_{ijjkk}).
\end{aligned}
\end{equation}
Similarly, by substituting both \(\xi_{jk}=-\xi_{kj}\) and \(\xi_{iijk}=-\xi_{kjii}\)  (reversal property) then by re-indexing, we obtain,
\begin{equation}\label{eq:ij_ijkk}
\begin{aligned}
\sum_{1\leq i,j,k\leq d}\Big (4(\xi_{ij},\xi_{ijkk})+4(\xi_{jk},\xi_{iijk})\Big )&=\sum_{1\leq i,j,k\leq d}\Big (4(\xi_{ij},\xi_{ijkk})+4(\xi_{kj},\xi_{kjii})\Big )\\
&=\sum_{1\leq i,j,k\leq d}8(\xi_{ij},\xi_{ijkk}).
\end{aligned}
\end{equation}
Finally, by using \(\xi_{iij}=\xi_{jii}\), \(\xi_{iik}=\xi_{kii}\) and \(\xi_{jjk}=\xi_{kjj}\) (reversal property) before re-indexing we obtain,
\begin{equation}\label{eq:iij_jjk}
\begin{aligned}
\sum_{1\leq i,j,k\leq d}\Big (2(\xi_{iij},\xi_{jkk})+2(\xi_{ijj},\xi_{ikk})+2(\xi_{iik},\xi_{jjk})\Big )&=\sum_{1\leq i,j,k\leq d}\Big (2(\xi_{jii},\xi_{jkk})+2(\xi_{ijj},\xi_{ikk})+2(\xi_{kii},\xi_{kjj})\Big )\\
&=\sum_{1\leq i,j,k\leq d}6(\xi_{ijj},\xi_{ikk}).
\end{aligned}
\end{equation}
Substituting \eqref{eq:i_ijjkk}, \eqref{eq:ij_ijkk}, \eqref{eq:iij_jjk} into \eqref{eq:i2j2k2_k2},
\begin{equation}\label{eq:k2_iijjkk}
\begin{aligned}
\sum_{1\leq i,j,k\leq d}\e^{\star 2}(\epsilon_i^2\epsilon_j^2\epsilon_k^2)=\sum_{1\leq i,j,k\leq d}\Big (4(\xi_i&,\xi_{ijjkk})+2(\xi_j,\xi_{iijkk})+8(\xi_{ij},\xi_{ijkk})\\
&+4(\xi_{ik},\xi_{ijjk})+6(\xi_{ijj},\xi_{ikk})+4(\xi_{ijk},\xi_{ijk})\Big ).
\end{aligned}
\end{equation}
Next, for \(k=4\),
\begin{equation}\label{eq:k4_iijjkk_full}
\begin{aligned}
\e^{\star 4}(\epsilon_i^2\epsilon_j^2\epsilon_k^2)=2(&\xi_i,\xi_i,\xi_j,\xi_{jkk})+2(\xi_i,\xi_j,\xi_j,\xi_{ikk})+2(\xi_i,\xi_i,\xi_k,\xi_{jjk})+2(\xi_i,\xi_k,\xi_k,\xi_{jjk})\\
&+2(\xi_j,\xi_j,\xi_k,\xi_{iik})+2(\xi_j,\xi_k,\xi_k,\xi_{iij})+8(\xi_i,\xi_j,\xi_k,\xi_{ijk})\\
&+2(\xi_i,\xi_i,\xi_{jk},\xi_{jk})+2(\xi_j,\xi_j,\xi_{ik},\xi_{ik})+2(\xi_k,\xi_k,\xi_{ij},\xi_{ij})\\
&+8(\xi_{i},\xi_{j},\xi_{ik},\xi_{jk})+8(\xi_{i},\xi_{k},\xi_{ij},\xi_{jk})+8(\xi_{j},\xi_{k},\xi_{ij},\xi_{ik}) .
\end{aligned}
\end{equation}
By substituting \(\xi_{jjk}=\xi_{kjj}, \xi_{iik}=\xi_{kii}, \xi_{iij}=\xi_{jii}\) (reversal property) before re-indexing,
\begin{equation}\label{eq:i_k_k_ijj}
\begin{aligned}
\sum_{1\leq i,j,k\leq d}12(\xi_i,\xi_k,\xi_k,\xi_{ijj})=\sum_{1\leq i,j,k\leq d}&\Big (2(\xi_i,\xi_i,\xi_j,\xi_{jkk})+2(\xi_i,\xi_j,\xi_j,\xi_{ikk})+2(\xi_i,\xi_i,\xi_k,\xi_{jjk})\\
&+2(\xi_i,\xi_k,\xi_k,\xi_{jjk})+2(\xi_j,\xi_j,\xi_k,\xi_{iik})+2(\xi_j,\xi_k,\xi_k,\xi_{iij})\Big ) .
\end{aligned}
\end{equation}
Also by only re-indexing,
\begin{equation}\label{eq:i_i_jk_jk}
\sum_{1\leq i,j,k\leq d}6(\xi_i,\xi_i,\xi_{jk},\xi_{jk})=\sum_{1\leq i,j,k\leq d}\Big (2(\xi_i,\xi_i,\xi_{jk},\xi_{jk})+2(\xi_j,\xi_j,\xi_{ik},\xi_{ik})+2(\xi_k,\xi_k,\xi_{ij},\xi_{ij})\Big ).
\end{equation}
Finally, by substituting \(\xi_{ij}=-\xi_{ji}\) (reversal / symmetric property) before re-indexing,
\begin{equation}\label{eq:i_j_ik_jk}
\begin{aligned}
\sum_{1\leq i,j,k\leq d}\Big ((\xi_{i},\xi_{k},\xi_{ij},\xi_{jk})+(\xi_{j},\xi_{k},\xi_{ij},\xi_{ik})\Big )&=\sum_{1\leq i,j,k\leq d}\Big (-(\xi_{i},\xi_{k},\xi_{ji},\xi_{jk})+(\xi_{j},\xi_{k},\xi_{ij},\xi_{ik})\Big )\\
&=0.
\end{aligned}
\end{equation}
We can eliminate one further term by using \(\xi_{ijk}+\xi_{jki}+\xi_{kij}=0\) (cylic property),
\begin{equation}\label{eq:i_j_k_ijk}
\begin{aligned}
\sum_{1\leq i,j,k\leq d}(\xi_i,\xi_j,\xi_k,\xi_{ijk})&=\frac 13 \sum_{1\leq i,j,k\leq d}\Big ((\xi_i,\xi_j,\xi_k,\xi_{ijk}) + (\xi_i,\xi_j,\xi_k,\xi_{ijk}) + (\xi_i,\xi_j,\xi_k,\xi_{ijk})\Big )\\
&=\frac 13 \sum_{1\leq i,j,k\leq d}\Big ((\xi_i,\xi_j,\xi_k,\xi_{ijk}) + (\xi_j,\xi_k,\xi_i,\xi_{jki}) + (\xi_k,\xi_i,\xi_j,\xi_{kij})\Big )\\
&=\frac 13 \sum_{1\leq i,j,k\leq d}(\xi_i,\xi_j,\xi_k,\xi_{ijk}+\xi_{jki}+\xi_{kij})\\
&=0.
\end{aligned}
\end{equation}
Substituting \eqref{eq:i_k_k_ijj}, \eqref{eq:i_i_jk_jk}, \eqref{eq:i_j_ik_jk} and \eqref{eq:i_j_k_ijk}, into \eqref{eq:k4_iijjkk_full} gives,
\begin{equation}\label{eq:k4_iijjkk}
\sum_{1\leq i,j,k\leq d}\e^{\star 4}(\epsilon_i^2\epsilon_j^2\epsilon_k^2)=\sum_{1\leq i,j,k\leq d}\Big (12(\xi_i,\xi_k,\xi_k,\xi_{ijj})+6(\xi_i,\xi_i,\xi_{jk},\xi_{jk})+8(\xi_{i},\xi_{j},\xi_{ik},\xi_{jk})\Big ) .
\end{equation}
For \(k=6\),
\begin{equation}\label{eq:k6_iijjkk}
\sum_{1\leq i,j,k\leq d}\e^{\star 6}(\epsilon_i^2\epsilon_j^2\epsilon_k^2)=\sum_{1\leq i,j,k\leq d}(\xi_i,\xi_i,\xi_j,\xi_j,\xi_k,\xi_k) .
\end{equation}
Substituting \eqref{eq:k2_iijjkk}, \eqref{eq:k4_iijjkk} and \eqref{eq:k6_iijjkk} into \eqref{eq:proof_ei2_ej2_ek2} we obtain the required result.
\end{proof}

\section{Eulerian idempotent in Lyndon basis}\label{app:lyndon}
In this appendix we state the Eulerian idempotents appearing in the cubature formulae constructed in \autoref{theorem:deg_7_cubature_formula} in the Lyndon basis. 
\begin{align*}
\e(\epsilon_i)&=\epsilon_i\\\\
\e(\epsilon_i\epsilon_j)&=\frac{1}{2}[\epsilon_i,\epsilon_j]\\\\
\e(\epsilon_i\epsilon_j\epsilon_k)&=\frac{1}{6}\Big (\big [\epsilon_i,[\epsilon_j,\epsilon_k]\big]+\big [[\epsilon_i,\epsilon_j],\epsilon_k\big] \Big )\\\\
\e(\epsilon_i\epsilon_j\epsilon_k\epsilon_l)&=\frac{1}{12}\Big (\big [\epsilon_i,\big [[\epsilon_j,\epsilon_k],\epsilon_l\big ]\big ]+\big [[\epsilon_i,\epsilon_j],[\epsilon_k,\epsilon_l]\big ]+\big [[\epsilon_i,\epsilon_k],[\epsilon_j,\epsilon_l]\big ]+\big [\big [\epsilon_i,[\epsilon_j,\epsilon_k]\big],\epsilon_l\big ]\Big )\\\\
\e(\epsilon_i\epsilon_j\epsilon_k\epsilon_l\epsilon_m)&=\frac{1}{60}\Big (12[\epsilon_i, [\epsilon_j, [\epsilon_k, [\epsilon_l, \epsilon_m]]]]+9[\epsilon_i, [\epsilon_j, [[\epsilon_k, \epsilon_m], \epsilon_l]]]+6[\epsilon_i, [[\epsilon_j, \epsilon_l], [\epsilon_k, \epsilon_m]]]\\
&+6[\epsilon_i, [[\epsilon_j, [\epsilon_l, \epsilon_m]], \epsilon_k]]+6[\epsilon_i, [[\epsilon_j, \epsilon_m], [\epsilon_k, \epsilon_l]]]+2[\epsilon_i, [[[\epsilon_j, \epsilon_m], \epsilon_l], \epsilon_k]]\\
&+3[[\epsilon_i, \epsilon_k], [\epsilon_j, [\epsilon_l, \epsilon_m]]]+[[\epsilon_i, \epsilon_k], [[\epsilon_j, \epsilon_m], \epsilon_l]]+2[[\epsilon_i, [\epsilon_k, \epsilon_l]], [\epsilon_j, \epsilon_m]]\\
&+3[[\epsilon_i, [\epsilon_k, [\epsilon_l, \epsilon_m]]], \epsilon_j]+2[[\epsilon_i, [\epsilon_k, \epsilon_m]], [\epsilon_j, \epsilon_l]]+[[\epsilon_i, [[\epsilon_k, \epsilon_m], \epsilon_l]], \epsilon_j]\\
&+3[[\epsilon_i, \epsilon_l], [\epsilon_j, [\epsilon_k, \epsilon_m]]]+[[\epsilon_i, \epsilon_l], [[\epsilon_j, \epsilon_m], \epsilon_k]]-[[[\epsilon_i, \epsilon_l], \epsilon_k], [\epsilon_j, \epsilon_m]]\\
&-[[[\epsilon_i, \epsilon_l], [\epsilon_k, \epsilon_m]], \epsilon_j]+2[[\epsilon_i, [\epsilon_l, \epsilon_m]], [\epsilon_j, \epsilon_k]]-[[[\epsilon_i, [\epsilon_l, \epsilon_m]], \epsilon_k], \epsilon_j]\\
&+3[[\epsilon_i, \epsilon_m], [\epsilon_j, [\epsilon_k, \epsilon_l]]]+[[\epsilon_i, \epsilon_m], [[\epsilon_j, \epsilon_l], \epsilon_k]]-[[[\epsilon_i, \epsilon_m], \epsilon_k], [\epsilon_j, \epsilon_l]]\\
&-[[[\epsilon_i, \epsilon_m], [\epsilon_k, \epsilon_l]], \epsilon_j]-[[[\epsilon_i, \epsilon_m], \epsilon_l], [\epsilon_j, \epsilon_k]]-2[[[[\epsilon_i, \epsilon_m], \epsilon_l], \epsilon_k], \epsilon_j]\Big )
\end{align*}

\section{A numerical toy example for linear SDEs}\label{app:numerics}

We present a simple numerical toy example demonstrating the application of our cubature formulas to linear stochastic differential equations.  This example serves to illustrate and confirm the convergence of our degree-7 cubature formula in comparison with lower-degree cubature approximations.  The convergence and practical implementation of cubature on Wiener space (up to degree-5 in general dimension) for more general test functions has been studied in Gyurko, Lyons \cite{gyurko}.

Specifically, we consider linear Stratonovich SDEs driven by three-dimensional noise, with a two-dimensional solution (a case that generalises readily), of the form
\[
\mathrm{d} Y^k = A^k_{i\gamma}Y^i \circ \mathrm{d}W^\gamma + B^k_iY^i \mathrm{d}t, \qquad Y_0 = y_0,
\]
in which we have used the Einstein summation convention. Its mean $\mathbb E Y$ can be computed by converting into It\^o form and passing to the expectation, yielding an ODE
\[
\frac{\dif \mathbb E Y^k}{\dif t} = \frac 12 \sum_{\gamma = 1}^d A^k_{i \gamma} A^i_{j \gamma} \mathbb E Y^j + B^k_i \mathbb E Y^i .
\]
\begin{wrapfigure}{r}{0.5\textwidth}
 \vspace{-20pt}
  \centering
  \includegraphics[width=0.5\textwidth]{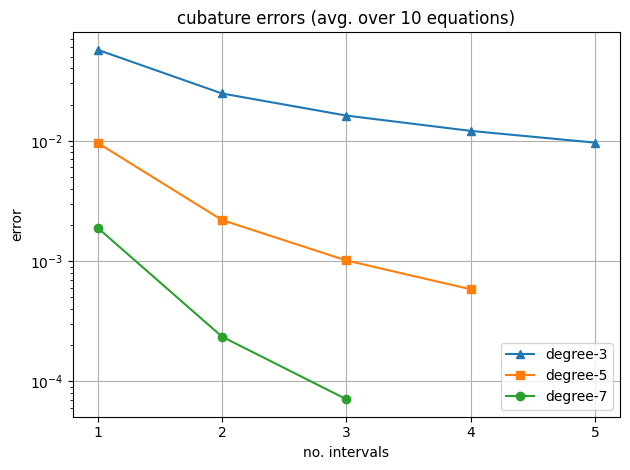}
  \vspace{-20pt}
  \caption{Cubature errors}
  \label{fig:plot}
   \vspace{-10pt}
\end{wrapfigure}
Solving this ODE to high accuracy provides us with an exact reference solution relative to which errors are calculated; it should be noted that this method for averaging is only available for linear SDEs and test functions. \par Our numerical experiments approximate $\mathbb E Y_1$ using degree-$3$ and $5$ formulas based on Lyons-Victoir cubature (\cite{LV04} and the degree-$7$ formula with support of size $648$ constructed in \autoref{sec:brokensym} with an Taylor scheme on uniform partitions of the time interval $[0,1]$. Relative errors $|\mathbb E Y^\mathrm{cub}_1 - \mathbb E Y^\mathrm{ODE}_1| / |\mathbb E Y^\mathrm{ODE}_1|$ are averaged across $10$ random choices of the triple $(y_0, A, B)$, with each entry normalised to have Euclidean norm $1$. Results are summarised in the plot \autoref{fig:plot} and the associated code may be found in the notebook \cite[\texttt{cubature\_plots.ipynb}]{code}.

The implementation of our toy example is naive, computing the full tree of the cubature approximation within memory constraints for three to five steps depending on the order of the cubature. In practice, exponential growth of the approximation tree can be overcome by combining cubature with partial sampling schemes such as TBBA \cite{CL02}, or recombination \cite{LL12,ninomiya} while preserving high-order accuracy; such an implementation lies beyond the scope of this appendix.

\paragraph{Acknowledgements}
EF was supported by
the EPSRC [EP/S026347/1]. TH was supported by the Crankstart Internship Programme. CL was supported by the Engineering and Physical Sciences Research Council (grant EP/V005413/1). TL was supported in part by UK Research and Innovation (UKRI) through the Engineering and Physical Sciences Research Council (EPSRC) via Programme Grants [Grant No. UKRI1010: High order mathematical and computational infrastructure for streamed data that enhance contemporary generative and large language models], [Grant No. EP/S026347/1: Unparameterised multi-model data, high order signatures and the mathematics of data science], and the UKRI AI for Science award [Grant No. UKRI2385: Creating Foundational Benchmarks for AI in Physical and Biological Complexity]. He was also supported by The Alan Turing Institute under the Defence and Security Programme (funded by the UK Government) and through the provision of research facilities; by the UK Government; and through CIMDA@Oxford, part of the AIR@InnoHK initiative funded by the Innovation and Technology Commission, HKSAR Government.

For the purpose of open access, the authors have applied a Creative Commons Attribution (CC BY) licence to any Author Accepted Manuscript version arising from this submission.

\bibliographystyle{alpha}
\bibliography{refs}

\end{document}